# FLUID LIMITS FOR NETWORKS WITH BANDWIDTH SHARING AND GENERAL DOCUMENT SIZE DISTRIBUTIONS

By H. Christian Gromoll[1] and Ruth J. Williams[2]

*University of Virginia and University of California, San Diego*

We consider a stochastic model of Internet congestion control, introduced by Massoulié and Roberts [*Telecommunication Systems* **15** (2000) 185–201], that represents the randomly varying number of flows in a network where bandwidth is shared among document transfers. In contrast to an earlier work by Kelly and Williams [*Ann. Appl. Probab.* **14** (2004) 1055–1083], the present paper allows interarrival times and document sizes to be generally distributed, rather than exponentially distributed. Furthermore, we allow a fairly general class of bandwidth sharing policies that includes the weighted $\alpha$-fair policies of Mo and Walrand [*IEEE/ACM Transactions on Networking* **8** (2000) 556–567], as well as certain other utility based scheduling policies. To describe the evolution of the system, measure valued processes are used to keep track of the residual document sizes of all flows through the network. We propose a fluid model (or formal functional law of large numbers approximation) associated with the stochastic flow level model. Under mild conditions, we show that the appropriately rescaled measure valued processes corresponding to a sequence of such models (with fixed network structure) are tight, and that any weak limit point of the sequence is almost surely a fluid model solution. For the special case of weighted $\alpha$-fair policies, we also characterize the invariant states of the fluid model.

**1. Introduction.** Massoulié and Roberts [18] have introduced and studied a model of Internet congestion control that represents the randomly

Received August 2006; revised March 2008.
[1]Supported in part by an NSF Mathematical Sciences Postdoctoral Research Fellowship, NSF Grant DMS-FRG-02-44323 and a European Union Marie Curie Postdoctoral Research Fellowship.
[2]Supported in part by NSF Grants DMS-03-05272 and DMS-06-04537.
*AMS 2000 subject classifications.* Primary 60K30; secondary 60F17, 90B15.
*Key words and phrases.* Bandwidth sharing, $\alpha$-fair, flow level Internet model, congestion control, simultaneous resource possession, fluid model, workload, measure valued process, invariant manifold.

varying number of flows in a network where bandwidth is shared dynamically among flows. The flows correspond to continuous transfers of individual elastic documents. This connection level model assumes a "separation of time scales" such that the time scale of the flow dynamics (of document arrivals and departures) is much longer than the time scale of the packet level dynamics on which rate control schemes such as TCP converge to equilibrium.

Subsequent to the work of Massoulié and Roberts [18], assuming exponentially distributed document sizes, de Veciana, Lee and Konstantopoulos [5] and Bonald and Massoulié [1] studied the stability of the flow level model operating under various bandwidth sharing policies. A bandwidth sharing policy generalizes the notion of a processor sharing discipline from a single resource to a network with several shared resources. Lyapunov functions constructed in [5] for weighted max-min fair and proportionally fair policies, and in [1] for weighted $\alpha$-fair policies ($\alpha \in (0, \infty)$) [19], imply positive recurrence of the Markov chain associated with the model when the average load on each resource is less than its capacity. Several authors [9, 13, 15, 16, 21, 22, 23] have considered variants of the Massoulié and Roberts model [18] and more general bandwidth sharing policies. In particular, Lin, Shroff and Srikant [15, 16, 21] have given sufficient conditions for stability where the assumption of time scale separation is relaxed. Ye [22], Ye, Ou and Yuan [23] and Hansen, Reynolds and Zachary [9] have given conditions for stability and instability with more general bandwidth sharing policies. Key and Massoulié [13] have considered a model with file transfers and streaming flows, certain utility based policies and relaxed capacity constraints. However, all of these works maintain a critical exponential distributional assumption on document sizes or holding times to enable the use of a relatively simple Markovian model. A major aim of our work is to relax this exponential assumption.

Here, we consider the model of Massoulié and Roberts, with generally distributed document sizes and interarrival times, operating under a fairly general bandwidth sharing policy. Important examples of this policy include the weighted $\alpha$-fair policies introduced by Mo and Walrand [19], and more generally certain utility based policies (see, e.g., [3, 13, 22, 23]) in the context of flow level models. We are interested in the stability and heavy traffic behavior of this flow level model. (Despite the claim in [1], the proof of sufficient conditions for stability under weighted $\alpha$-fair policies given there does not apply when document sizes are other than exponentially distributed. The reason for this is that the method of Dai [4] quoted there implicitly assumes (through the form of the model equations) that the service discipline is a head-of-the-line discipline. Consequently, the method does not apply in general to processor sharing type disciplines, such as the bandwidth sharing policies considered here. In the case of exponentially distributed document sizes, one can equate the distribution of the queue length process for a bandwidth sharing model with the queue length process of a stochastic processing network (cf. [10]) operating under a head-of-the-line policy. Even then, to conclude the stability result using an analogue of Dai's result, one has to generalize the results of [4] to stochastic processing networks from multiclass queueing networks. However, in the case of exponential interarrival times and document sizes, the Lyapunov function given in [1] can be used directly on the original Markov chain stochastic model to establish stability under the nominal condition that the average load placed on each resource is less than its capacity.)

There are a few results on sufficient conditions for stability of the flow level model with general document size distributions. With Poisson arrivals and document sizes having a phase-type distribution, for a weighted $\alpha$-fair policy with $\alpha = 1$, Lakshmikantha, Beck and Srikant [14] have established stability of some two resource linear networks and a $2 \times 2$ grid network when the average load on each resource is less than its capacity. For generally distributed interarrival and document sizes, Bramson [2] has shown sufficiency of such a condition for stability under a max-min fair policy (corresponding to an $\alpha$-fair policy as $\alpha \to \infty$). Under proportional fair sharing, Massoulié [17] has recently established stability of a fluid model for the flow level model with exponential interarrival and document sizes, and additional routing. From this, he infers stability of the stochastic flow level model when documents have phase-type distributions. In general, however, it remains an open question whether, with renewal arrivals and arbitrarily (rather than exponentially) distributed document sizes, the flow level model is stable under a weighted $\alpha$-fair (or more general) bandwidth sharing policy when the nominal load placed on each resource is less than its capacity. In contemporaneous work, Chiang, Shah and Tang [3] have developed a fluid approximation for the flow level model when the arrival rate and capacity are allowed to grow proportionally but the bandwidth per flow stays uniformly bounded. Using their fluid model, they derive some conclusions concerning rate stability for general (bounded) document size distributions when $\alpha \in (0, \infty)$ is sufficiently small.

This paper is a first step in our study of the flow level model with general interarrival and document size distributions, and a general bandwidth sharing policy. Here, we define measure valued processes that keep track of the residual sizes of all documents in the system at any given time. We propose a fluid model (or formal functional law of large numbers approximation) associated with the stochastic flow level model. Under mild conditions, we show that the measure valued processes corresponding to a fluid scaled sequence of such models (with fixed network structure) are tight and that any weak limit point of the sequence is almost surely a fluid model solution. For weighted $\alpha$-fair policies, we also characterize the invariant states for the fluid model. In future work, we plan to study the asymptotic behavior of fluid model solutions and to use that to study the stability and heavy traffic behavior of the associated flow level models. A summary of the results of this paper as they pertain to weighted $\alpha$-fair policies appears in [8], along with two examples showing stability of the fluid model under a natural condition for linear networks and simple tree networks.

The paper is organized as follows. In Section 2, we define the network structure, the bandwidth sharing policy, the stochastic flow level model and we introduce the measure valued processes used to describe the evolution of the system. The notion of a fluid model solution is defined in Section 3. In Section 4, we introduce a sequence of flow level models and state our main result concerning the tightness of this sequence and that weak limit points are fluid model solutions (see Theorem 4.1). The proof of the main result is given in Section 5. In Section 6, we characterize the invariant states of the fluid model for weighted $\alpha$-fair policies.

1.1. *Notation.* Let $\mathbb{N} = \{1, 2, \ldots\}$, let $\mathbb{R} = (-\infty, \infty)$, and let $\mathbb{R}^d$ denote $d$-dimensional Euclidean space. For $x, y \in \mathbb{R}$, $x \vee y$ is the maximum of $x$ and $y$, $x \wedge y$ is the minimum of $x$ and $y$, $x^+$ is the positive part and $\lfloor x \rfloor$ is the integer part of $x$. For $x, y \in \mathbb{R}^d$, let $\|x\| = \max_{i=1}^d |x_i|$, and interpret vector inequalities componentwise: $x \leq y$ means $x_i \leq y_i$ for all $i = 1, \ldots, d$. The positive $d$-dimensional orthant is denoted $\mathbb{R}^d_+ = \{x \in \mathbb{R}^d : x \geq 0\}$. To ease notation throughout the paper, all vectors are considered to be column vectors when used in mathematical expressions, but will be written out as row vectors within paragraphs. Also, define $c/0$ to be zero for any real constant $c$, and define a sum over an empty set of indices or of the form $\sum_{k=j}^{l}$ with $j > l$ to be zero.

For two functions $f$ and $g$ with the same domain, $f \equiv g$ means $f(x) = g(x)$ for all $x$ in the domain. For a bounded function $f : \mathbb{R}_+ \to \mathbb{R}$, let $\|f\|_\infty = \sup_{x \in \mathbb{R}_+} |f(x)|$. Let $\mathbf{C}_b(\mathbb{R}_+)$ be the set of continuous bounded functions $f : \mathbb{R}_+ \to \mathbb{R}$, let $\mathbf{C}^1(\mathbb{R}_+)$ be the set of once continuously differentiable functions $f : \mathbb{R}_+ \to \mathbb{R}$, and let $\mathbf{C}^1_b(\mathbb{R}_+)$ be the set of functions $f$ in $\mathbf{C}^1(\mathbb{R}_+)$ that together with the first derivative $f'$, are bounded on $\mathbb{R}_+$. If $w \in \mathbf{C}^1_b(\mathbb{R}_+)$ is a function of time, its derivative will be denoted by $\dot{w}$. For a Polish (complete separable metric) space $\mathcal{S}$, let $\mathbf{D}([0, \infty), \mathcal{S})$ be the space of right continuous functions from $[0, \infty)$ into $\mathcal{S}$ that have left limits in $\mathcal{S}$. Endow this space with the Skorohod $J_1$-topology. For a finite nonnegative Borel measure $\xi$ on $\mathbb{R}_+$ and a $\xi$-integrable function $f : \mathbb{R}_+ \to \mathbb{R}$, define

$$\langle f, \xi \rangle = \int_{\mathbb{R}_+} f \, d\xi.$$

If $\xi = (\xi_1, \ldots, \xi_d)$ is a vector of such measures, then $\langle f, \xi \rangle$ is the vector $(\langle f, \xi_1 \rangle, \ldots, \langle f, \xi_d \rangle)$. All functions $f : \mathbb{R}_+ \to \mathbb{R}$ are extended to be identically zero on $(-\infty, 0)$ so that $f(\cdot - x)$ is well defined on $\mathbb{R}_+$ for all $x > 0$. Let $\chi : \mathbb{R}_+ \to \mathbb{R}_+$ denote the identity function $\chi(x) = x$ for $x \in \mathbb{R}_+$.

Let $\mathbf{M}$ be the set of finite nonnegative Borel measures on $\mathbb{R}_+$, endowed with the weak topology: $\xi^k \xrightarrow{\mathbf{w}} \xi$ in $\mathbf{M}$ if and only if $\langle f, \xi^k \rangle \to \langle f, \xi \rangle$ for all $f \in \mathbf{C}_b(\mathbb{R}_+)$. This topology is induced by the following generalization of the Prohorov metric: for $\xi, \zeta \in \mathbf{M}$ define

$$\begin{aligned}(1.1)\quad \mathbf{d}[\xi, \zeta] = \inf\{\varepsilon > 0 : {}&\xi(B) \leq \zeta(B^\varepsilon) + \varepsilon \text{ and} \\ &\zeta(B) \leq \xi(B^\varepsilon) + \varepsilon \text{ for all nonempty closed } B \subset \mathbb{R}_+\},\end{aligned}$$

where $B^\varepsilon = \{x \in \mathbb{R}_+ : \inf_{y \in B} |x - y| < \varepsilon\}$. It will be convenient to extend the notion of uniform integrability for random variables (and their associated distributions) to elements of $\mathbf{M}$. Call a sequence $\{\xi^k\} \subset \mathbf{M}$ uniformly integrable, if $\langle \chi, \xi^k \rangle < \infty$ for all $k$ and

$$\lim_{x \to \infty} \sup_k \langle \chi 1_{[x, \infty)}, \xi^k \rangle = 0.$$

It is easy to show that if $\{\xi^k\} \subset \mathbf{M}$ is uniformly integrable and $\xi^k \xrightarrow{\mathbf{w}} \xi$, then $\langle \chi, \xi \rangle < \infty$ and $\langle \chi, \xi^k \rangle \to \langle \chi, \xi \rangle$.

For $\mathbf{I} \in \mathbb{N}$, let

$$\mathbf{M}^\mathbf{I} = \{(\xi_1, \ldots, \xi_\mathbf{I}) : \xi_i \in \mathbf{M} \text{ for all } i \leq \mathbf{I}\}$$

and for $\xi, \zeta \in \mathbf{M}^\mathbf{I}$, define

$$(1.2)\qquad \mathbf{d}_\mathbf{I}[\xi, \zeta] = \max_{i \leq \mathbf{I}} \mathbf{d}[\xi_i, \zeta_i].$$

Equipped with the metric $\mathbf{d}_\mathbf{I}[\cdot, \cdot]$, the space $\mathbf{M}^\mathbf{I}$ is Polish. Convergence of a sequence $\{\xi^k\}$ to $\xi$ in $\mathbf{M}^\mathbf{I}$ is also denoted $\xi^k \xrightarrow{\mathbf{w}} \xi$. The zero measure in $\mathbf{M}$ is denoted by $\mathbf{0}$.

The notation $X \sim Y$ means $X$ and $Y$ are equal in distribution, and $X^n \Rightarrow X$ means the sequence $\{X^n\}$ converges in distribution to $X$. All continuous time stochastic processes used in this work are assumed to have sample paths that are right continuous with left limits.

**2. Flow level model.** This section defines the network structure, the bandwidth sharing policy and the stochastic flow level model.

2.1. *Network structure.* Consider a network with finitely many *resources* labelled by $j = 1, \ldots, \mathbf{J}$, and a finite set of *routes* labeled by $i = 1, \ldots, \mathbf{I}$. A route $i$ is a nonempty subset of $\{1, \ldots, \mathbf{J}\}$, interpreted as the set of resources used by the route. Let $A$ be the $\mathbf{J} \times \mathbf{I}$ incidence matrix satisfying $A_{ji} = 1$ if resource $j$ is used by route $i$, and $A_{ji} = 0$ otherwise. Since each route is a nonempty subset of $\{1, \ldots, \mathbf{J}\}$, no column of $A$ is identically zero.

A *flow* on route $i$ is the continuous transfer of a document through the resources used by the route. Assume that while being transferred, a flow takes simultaneous possession of all resources on its route. The *processing rate* allocated to a flow is the rate at which the associated document is being transferred. There may be multiple flows on a route, and the *bandwidth* $\Lambda_i$ allocated to route $i$ is the sum of the processing rates allocated to flows on route $i$. The *bandwidth allocated through resource* $j$ is the sum of the bandwidths allocated to routes using resource $j$. Assume that each resource $j \leq \mathbf{J}$ has finite *capacity* $C_j > 0$, interpreted as the maximum bandwidth that can be allocated through it. Let $C = (C_1, \ldots, C_\mathbf{J})$ be the vector of capacities in $\mathbb{R}_+^\mathbf{J}$. Then any vector $\Lambda = (\Lambda_1, \ldots, \Lambda_\mathbf{I})$ of bandwidth allocations must satisfy

$$A\Lambda \leq C.$$

2.2. *Bandwidth sharing policy.* We consider the network operating under a policy that dynamically allocates bandwidth to routes as a function of the number of flows on all routes. The resulting allocation to each route is shared equally among individual flows on that route.

Let $Z_i(t)$ denote the number of flows on route $i \leq \mathbf{I}$ at time $t$, and let $Z(t) = (Z_1(t), \ldots, Z_\mathbf{I}(t))$ be the corresponding vector in $\mathbb{R}_+^\mathbf{I}$. The bandwidth allocated to route $i$ at time $t$ is a function of the vector $Z(t)$ and is denoted by $\Lambda_i(Z(t))$. The corresponding vector of bandwidth allocations at time $t$ is $\Lambda(Z(t)) = (\Lambda_1(Z(t)), \ldots, \Lambda_\mathbf{I}(Z(t)))$. Although the coordinates of $Z(\cdot)$ are nonnegative and integer valued, we assume that the function $\Lambda$ is defined on the entire orthant $\mathbb{R}_+^\mathbf{I}$ to accommodate fluid analogues of $Z(\cdot)$ later.

DEFINITION 2.1. A bandwidth sharing policy for the network $(A, C)$ is a function $\Lambda : \mathbb{R}_+^\mathbf{I} \to \mathbb{R}_+^\mathbf{I}$ such that for each $z \in \mathbb{R}_+^\mathbf{I}$:

(i) $\Lambda_i(z) > 0$ for each $i$ such that $z_i > 0$,
(ii) $\Lambda_i(z) = 0$ for each $i$ such that $z_i = 0$,
(iii) $A\Lambda(z) \leq C$,
(iv) $\Lambda(rz) = \Lambda(z)$ for each $r > 0$, and such that for each $i \leq \mathbf{I}$,
(v) $\Lambda_i(\cdot)$ is continuous on $\{z \in \mathbb{R}_+^\mathbf{I} : z_i > 0\}$.

Properties (i) and (ii) imply that routes with active flows may not idle, and that no bandwidth is allocated to routes with no flows. Property (iii) is the basic feasibility constraint, and property (iv) requires that bandwidth allocations are invariant under scaling. Note that by property (iii), since each route uses at least one resource, we have

(2.1) $$\sup_{z \in \mathbb{R}_+^\mathbf{I}} \|\Lambda(z)\| \leq \|C\|.$$

We assume further that the bandwidth $\Lambda_i(Z(t))$ allocated to route $i$ at time $t$ is shared equally by all flows on the route. That is, if there are $Z_i(t) > 0$ flows on route $i$ at time $t$, then each flow is allocated a processing rate of $\Lambda_i(Z(t))/Z_i(t)$ at time $t$.

The following property of $\Lambda(\cdot)$ will be used later in this paper.

LEMMA 2.2. *Let $\Lambda(\cdot)$ be a bandwidth sharing policy for the network $(A, C)$. For each $\varepsilon, M \in (0, \infty)$, there exists $c > 0$ such that for each $i \leq \mathbf{I}$,*

$$\Lambda_i(z) \geq c \qquad \text{on } \{z \in \mathbb{R}_+^\mathbf{I} : z_i \geq \varepsilon, \|z\| \leq M\}.$$

PROOF. For each $i \leq \mathbf{I}$, the function $\Lambda_i(\cdot)$ is continuous and strictly positive on $\{z \in \mathbb{R}_+^\mathbf{I} : z_i > 0\}$ by Definition 2.1. So $\Lambda(\cdot)$ is bounded away from zero on the compact subset $\{z \in \mathbb{R}_+^\mathbf{I} : z_i \geq \varepsilon, \|z\| \leq M\}$. □

An important class of bandwidth sharing policies satisfying Definition 2.1 is described below.

EXAMPLE. The following family of policies was introduced by Mo and Walrand [19]. Fix a parameter $\alpha \in (0, \infty)$ and a vector of strictly positive weights $\kappa = (\kappa_1, \ldots, \kappa_\mathbf{I})$. For $z \in \mathbb{R}_+^\mathbf{I}$, let $\mathcal{I}_0(z) = \{i \leq \mathbf{I} : z_i = 0\}$ and $\mathcal{I}_+(z) = \{i \leq \mathbf{I} : z_i > 0\}$. Let $\mathbb{O}(z) = \{\lambda \in \mathbb{R}_+^\mathbf{I} : \lambda_i = 0 \text{ for all } i \in \mathcal{I}_0(z)\}$. Define a function $G_z : \mathbb{R}_+^\mathbf{I} \to [-\infty, \infty)$ by

$$(2.2) \qquad G_z(\lambda) = \begin{cases} \displaystyle\sum_{i \in \mathcal{I}_+(z)} \kappa_i z_i^\alpha \frac{\lambda_i^{1-\alpha}}{1-\alpha}, & \alpha \in (0, \infty) \setminus \{1\}, \\ \displaystyle\sum_{i \in \mathcal{I}_+(z)} \kappa_i z_i \log \lambda_i, & \alpha = 1, \end{cases}$$

where the value of $G_z(\lambda)$ is taken to be $-\infty$ if $\alpha \in [1, \infty)$ and $\lambda_i = 0$ for some $i \in \mathcal{I}_+(z)$, and $G_z(\lambda) = 0$ if $\mathcal{I}_+(z) = \varnothing$. For each $z \in \mathbb{R}_+^\mathbf{I}$, define $\Lambda(z)$ as the unique vector $\lambda \in \mathbb{R}_+^\mathbf{I}$ that solves the optimization problem:

(2.3) $\qquad\qquad\qquad$ maximize $G_z(\lambda)$,

(2.4) $\qquad\qquad\qquad$ subject to $A\lambda \leq C$,

(2.5) $\qquad\qquad\qquad$ over $\mathbb{O}(z)$.

The resulting allocation is called a *weighted $\alpha$-fair allocation*, and the function $\Lambda : \mathbb{R}_+^\mathbf{I} \to \mathbb{R}_+^\mathbf{I}$ is called a *weighted $\alpha$-fair bandwidth sharing policy*. Note that by (2.4) and (2.5), $\Lambda$ satisfies properties (ii) and (iii) of Definition 2.1. Properties (i), (iv), and (v) hold for $\Lambda$ by the proofs of Lemmas A.1–A.3 of [12]. (Although it is assumed at the beginning of [12] that $A$ has full row rank, scrutiny of the proofs of Lemmas A.1–A.3 in [12] reveals that this assumption is not used in verifying these properties.) When $\kappa_i = 1$ for all $i \leq \mathbf{I}$, the case $\alpha = 1$ and the limiting cases $\alpha \to 0$ and $\alpha \to \infty$ correspond, respectively, to a bandwidth allocation that is *proportionally fair*, achieves *maximum throughput*, or is *max-min fair* [1, 19].

Some authors (see, e.g., Ye [22], Ye, Ou and Yuan [23], Key and Massoulié [13] and Chiang, Shah and Tang [3]) have proposed more general objective functions than $G_z(\cdot)$ for determining bandwidth allocations in the context of flow level models. Indeed, the optimization problem (2.3)–(2.5) can be replaced by an equivalent one for the *per flow* bandwidth allocations $x_i = \lambda_i / z_i$ for $i \in \mathcal{I}_+(z)$, where $G_z(\lambda)$ given by (2.2) is replaced by

$$\sum_{i \in \mathcal{I}_+(z)} \kappa_i z_i \mathcal{U}(x_i)$$

and the utility function $\mathcal{U}$ is given by

$$\mathcal{U}(x) = \begin{cases} \dfrac{x^{1-\alpha}}{1-\alpha}, & \alpha \in (0, \infty) \setminus \{1\}, \\ \log(x), & \alpha = 1. \end{cases}$$

When a more general strictly concave utility function $\mathcal{U}$ is used, properties (ii) and (iii) are immediate from the form of the optimization problem, properties (i) and (v) will hold under suitable regularity conditions on $\mathcal{U}$, and (as pointed out by Chiang, Shah and Tang [3]), the critical scaling property (iv) will be satisfied if $\mathcal{U}$ has the scaling property that $\mathcal{U}(rx) = g(r)\mathcal{U}(x)$ for all $r > 0$, $x > 0$, and some function $g: (0, \infty) \to (0, \infty)$. As Chiang, Shah and Tang [3] also indicate, by seeking a scaling limit involving large capacities, one can relax this last assumption. However, this involves allowing the network capacity $C$ to grow with the scaling limit and is a different limiting regime than the one considered here; the present analysis is oriented toward a system with fixed network parameters $A, C$.

2.3. *Stochastic model.* Henceforth, we fix a network structure $(A, C)$ and a bandwidth sharing policy $\Lambda$. Our stochastic model of document flows consists of the following: a collection of stochastic primitives $E_1, \ldots, E_\mathbf{I}$ and $\{v_{1k}\}_{k=1}^{\infty}, \ldots, \{v_{\mathbf{I}k}\}_{k=1}^{\infty}$ describing the arrivals of document flows (including their sizes) to the network, a random initial condition $\mathcal{Z}(0) \in \mathbf{M}^\mathbf{I}$ specifying the state of the system at time zero and a collection of performance processes describing the time evolution of the system state. The performance processes are defined in terms of the primitives and initial condition through a set of descriptive equations. The random objects involved are defined on a common probability space $(\Omega, \mathscr{F}, \mathbf{P})$, with expectation operator $\mathbf{E}$.

The stochastic primitives consist of an *exogenous arrival process* $E_i$ and a sequence of *document sizes* $\{v_{ik}\}_{k=1}^{\infty}$ for each route $i \leq \mathbf{I}$. The arrival process $E_i$ is a counting process, that is, a nondecreasing, nonnegative integer valued process starting from zero. For $t \geq 0$, $E_i(t)$ represents the number of flows that have arrived to route $i$ during the time interval $(0, t]$. The $k$th such arrival is called flow $k$ on route $i$ and arrives at time $U_{ik} = \inf\{t \geq 0 : E_i(t) \geq k\}$ (note that simultaneous arrivals are allowed). Flows already on route $i$ at time zero are called *initial flows*.

For each $i \leq \mathbf{I}$ and $k \geq 1$, the random variable $v_{ik}$ represents the initial size of the document associated with flow $k$ on route $i$. This is the cumulative amount of processing that must be allocated to the flow to complete its transfer through the network. Assume that for each $i \leq \mathbf{I}$, the random variables $\{v_{ik}\}_{k=1}^{\infty}$ are strictly positive and form a sequence of independent and identically distributed random variables with common distribution $\vartheta_i$ on $\mathbb{R}_+$. Assume that the mean $\langle \chi, \vartheta_i \rangle \in (0, \infty)$ and let $\mu_i = \langle \chi, \vartheta_i \rangle^{-1}$. We make no further assumptions about the relationship between $\mu_i$ and $E_i$. The fluid approximation result stated in Section 4.3 below is valid for both underloaded and overloaded systems.

It will be convenient to combine the collection of stochastic primitives into a single, measure valued load process. For each $x \in \mathbb{R}_+$, let $\delta_x \in \mathbf{M}$ denote the Dirac point measure at $x$.

DEFINITION 2.3. For $i \leq \mathbf{I}$, define the load process for route $i$ by

$$\mathcal{L}_i(t) = \sum_{k=1}^{E_i(t)} \delta_{v_{ik}}, \qquad t \geq 0. \tag{2.6}$$

For $t \geq s \geq 0$, define the increment $\mathcal{L}_i(s, t) = \mathcal{L}_i(t) - \mathcal{L}_i(s)$.

The process $\mathcal{L} = (\mathcal{L}_1, \ldots, \mathcal{L}_\mathbf{I})$ is a random element of the Skorohod space $\mathbf{D}([0, \infty), \mathbf{M}^\mathbf{I})$. Note that $\mathcal{L}(s, t) \in \mathbf{M}^\mathbf{I}$ for all $t \geq s \geq 0$.

The initial condition specifies $Z(0) = (Z_1(0), \ldots, Z_\mathbf{I}(0))$, the number of initial flows on each route at time zero, as well as the initial sizes of the documents associated to these flows. Assume that the components of $Z(0)$ are nonnegative, integer valued random variables. The initial document sizes of the initial flows on route $i \leq \mathbf{I}$ are the first $Z_i(0)$ elements of a sequence $\{\tilde{v}_{il}\}_{l=1}^{\infty}$ of strictly positive random variables. A convenient way to express the initial condition is to define an initial random vector of measures $\mathcal{Z}(0) \in \mathbf{M}^\mathbf{I}$ with components

$$\mathcal{Z}_i(0) = \sum_{l=1}^{Z_i(0)} \delta_{\tilde{v}_{il}}, \qquad i \leq \mathbf{I}.$$

Henceforth, $\mathcal{Z}(0)$ will be used as the initial condition for the network.

The performance processes consist of a measure valued process $\mathcal{Z}$, taking values in $\mathbf{D}([0, \infty), \mathbf{M}^\mathbf{I})$, and a collection of auxiliary processes $(Z, T, U, W)$.

The process $Z = (Z_1, \ldots, Z_\mathbf{I})$ takes values in $\mathbf{D}([0,\infty), \mathbb{R}_+^\mathbf{I})$. For $i \leq \mathbf{I}$ and $t \geq 0$, $Z_i(t)$ is the number of (active) flows on route $i$ at time $t$. Recall that at time $t$, the bandwidth allocated to route $i$ is $\Lambda_i(Z(t))$, and this bandwidth is shared equally by all $Z_i(t)$ flows on route $i$; each such flow receives a processing rate of $\Lambda_i(Z(t))/Z_i(t)$, which equals zero by convention if $Z_i(t) = 0$. Thus, a flow that is active on route $i$ during a time interval $[s,t] \subset [0,\infty)$ receives *cumulative service during* $[s,t]$ equal to

$$(2.7) \qquad S_i(s,t) = \int_s^t \frac{\Lambda_i(Z(u))}{Z_i(u)} \, du.$$

Consider flow $k$ on route $i$. This flow arrives at time $U_{ik}$ and has initial document size $v_{ik}$. At time $t \geq U_{ik}$, the cumulative service received by this flow during $[U_{ik}, t]$ equals $S_i(U_{ik}, t) \wedge v_{ik}$. The amount of service still required therefore equals $(v_{ik} - S_i(U_{ik}, t))^+$. (Once this latter quantity becomes zero, the flow becomes inactive, i.e., it departs from the system.) A similar description applies for the initial flows on route $i$. For $t \geq 0$, $k \leq E_i(t)$, and $l \leq Z_i(0)$, define the *residual document size at time* $t$ of flow $k$ on route $i$ and initial flow $l$ on route $i$, by

$$(2.8) \qquad v_{ik}(t) = (v_{ik} - S_i(U_{ik}, t))^+ \quad \text{and} \quad \tilde{v}_{il}(t) = (\tilde{v}_{il} - S_i(0,t))^+,$$

respectively.

The measure valued process $\mathcal{Z} = (\mathcal{Z}_1, \ldots, \mathcal{Z}_\mathbf{I})$ is called the *state descriptor*; it tracks the residual document sizes of flows on all routes at any given time. Let $\delta_x^+ \in \mathbf{M}$ denote the Dirac measure at $x$ if $x \in (0,\infty)$, with $\delta_0^+ = \mathbf{0}$. For $t \geq 0$ and $i \leq \mathbf{I}$, define the finite Borel measure

$$(2.9) \qquad \mathcal{Z}_i(t) = \sum_{l=1}^{Z_i(0)} \delta_{\tilde{v}_{il}(t)}^+ + \sum_{k=1}^{E_i(t)} \delta_{v_{ik}(t)}^+.$$

Note that at $t = 0$, this definition coincides with the definition of the initial condition $\mathcal{Z}(0)$. Note also that by definition of the residual document sizes, the measure $\mathcal{Z}_i(t)$ has a unit of mass only for flows on route $i$ that have not yet completed transfer. Thus, for all $t \geq 0$ and $i \leq \mathbf{I}$,

$$(2.10) \qquad Z_i(t) = \langle 1, \mathcal{Z}_i(t) \rangle.$$

For $t \geq 0$ and $i \leq \mathbf{I}$, define

$$(2.11) \qquad T_i(t) = \int_0^t \Lambda_i(Z(s)) \, ds.$$

The process $T$ takes values in $\mathbf{D}([0,\infty), \mathbb{R}_+^\mathbf{I})$ and tracks the cumulative bandwidth allocated to each route. For $t \geq 0$, define

$$(2.12) \qquad U(t) = Ct - AT(t).$$

The process $U$ takes values in $\mathbf{D}([0,\infty), \mathbb{R}_+^{\mathbf{J}})$ and tracks the cumulative unused bandwidth capacity of each resource. Since $A\Lambda(z) \leq C$ for all $z \in \mathbb{R}_+^{\mathbf{I}}$, the process $U$ is nondecreasing. For $t \geq 0$, define

$$(2.13) \qquad W(t) = \langle \chi, \mathcal{Z}(t) \rangle.$$

Recall that $\chi(x) = x$ and that integration against the vector of measures $\mathcal{Z}(t)$ is interpreted componentwise. The process $W$ takes values in the path space $\mathbf{D}([0,\infty), \mathbb{R}_+^{\mathbf{I}})$. By (2.9), $W_i(t)$ is the sum of all residual document sizes on route $i$ at time $t$. Thus, $W_i(t)$ represents the immediate amount of work still to be transferred on route $i$ at time $t$. It can be shown that

$$(2.14) \qquad W_i(t) = W_i(0) + \langle \chi, \mathcal{L}_i(t) \rangle - T_i(t), \qquad i \leq \mathbf{I},\ t \geq 0.$$

This equation describes the workload on route $i$ at time $t$ in terms of the cumulative amount of work that arrives to and is processed on the route during $[0, t]$.

**3. Fluid model.** In this section, we define a fluid analogue of the stochastic model introduced in Section 2.3. The main goal of the paper is to establish, under mild assumptions, that a sequence of fluid scaled stochastic state descriptors is tight and that weak limit points are fluid model solutions (see Theorem 4.1 below). Fix a vector of strictly positive constants $\nu = (\nu_1, \ldots, \nu_{\mathbf{I}})$ and a vector of probability measures $\vartheta = (\vartheta_1, \ldots, \vartheta_{\mathbf{I}})$ in $\mathbf{M}^{\mathbf{I}}$, satisfying $\langle \chi, \vartheta_i \rangle < \infty$ and $\langle 1_{\{0\}}, \vartheta_i \rangle = 0$ for all $i \leq \mathbf{I}$. The constant $\nu_i$, $i \leq \mathbf{I}$, will be the fluid analogue of mean arrival rate to route $i$ in the stochastic model (when that exists). Let $\mu_i = \langle \chi, \vartheta_i \rangle^{-1}$ and $\rho_i = \nu_i / \mu_i$ for each $i \leq \mathbf{I}$. We do not impose criticality assumptions on the constants $\rho_i$; they may take any value in $(0, \infty)$. The fluid model consists of a deterministic measure valued function of time, called the *fluid model solution*, and a collection of auxiliary functions of time defined below.

DEFINITION 3.1. Given a continuous function $\zeta : [0, \infty) \to \mathbf{M}^{\mathbf{I}}$, define the auxiliary functions $(z, \tau, u, w)$ of $\zeta$, with respect to the data $(A, C, \Lambda, \nu, \vartheta)$, by

$$z(t) = \langle 1, \zeta(t) \rangle,$$
$$\tau_i(t) = \int_0^t (\Lambda_i(z(s)) 1_{(0,\infty)}(z_i(s)) + \rho_i 1_{\{0\}}(z_i(s)))\, ds, \qquad i \leq \mathbf{I},$$
$$u(t) = Ct - A\tau(t),$$
$$w(t) = \langle \chi, \zeta(t) \rangle$$

for all $t \geq 0$.

Here $z(t)$ and $\tau(t)$ take values in $\mathbb{R}_+^{\mathbf{I}}$ and $u(t)$ will take values in $\mathbb{R}_+^{\mathbf{J}}$. On the other hand, $w(t)$ takes values in $[0, \infty]^{\mathbf{I}}$, as $\zeta(t)$ need not have a finite first moment [see (ii) below].

A fluid model solution is now defined via projections against test functions in the class
$$\mathcal{C} = \{f \in \mathbf{C}_b^1(\mathbb{R}_+) : f(0) = f'(0) = 0\}.$$

DEFINITION 3.2. *A fluid model solution for the data $(A, C, \Lambda, \nu, \vartheta)$ is a continuous function $\zeta : [0, \infty) \to \mathbf{M}^{\mathbf{I}}$ that, together with its auxiliary functions $(z, \tau, u)$, satisfies:*

(i) $\|\langle 1_{\{0\}}, \zeta(t)\rangle\| = 0$ for all $t \geq 0$,
(ii) $u_j$ is nondecreasing for all $j \leq \mathbf{J}$,
(iii) for each $f \in \mathcal{C}$, $i \leq \mathbf{I}$, and $t \geq 0$,

$$\langle f, \zeta_i(t)\rangle = \langle f, \zeta_i(0)\rangle - \int_0^t \langle f', \zeta_i(s)\rangle \frac{\Lambda_i(z(s))}{z_i(s)}\, ds$$
(3.1)
$$+ \nu_i \langle f, \vartheta_i \rangle \int_0^t 1_{(0,\infty)}(z_i(s))\, ds.$$

Recall that in (3.1), the integrand in the first integral term is defined to be zero when its denominator is zero.

In Definition 3.2, it is possible to extend property (iii) to the class of functions $\{f \in \mathbf{C}_b^1(\mathbb{R}_+) : f(0) = 0\}$, yielding an equivalent definition. The more restrictive class $\mathcal{C}$ is used here to facilitate parts of the proof of Theorem 4.1 below. In particular, since $f(0) = f'(0) = 0$ on $\mathcal{C}$, a function in $\mathcal{C}$ can be extended to a function in $\mathbf{C}_b^1(\mathbb{R})$ by defining it to be identically zero on $(-\infty, 0)$.

When the initial fluid workload is finite, we have the following result.

LEMMA 3.3. *Suppose $\zeta$ is a fluid model solution with finite initial workload, that is, $w_i(0) = \langle \chi, \zeta_i(0)\rangle < \infty$ for all $i \leq \mathbf{I}$. Then the fluid workload function $w$ associated with $\zeta$ satisfies the following for each $i \leq \mathbf{I}$ and $t \geq 0$:*

$$w_i(t) = w_i(0) + \int_0^t (\rho_i - \Lambda_i(z(s))) 1_{(0,\infty)}(z_i(s))\, ds$$
(3.2)
$$= w_i(0) + \rho_i t - \tau_i(t).$$

*In particular, the fluid workload $w_i(t)$ is finite for all $t \geq 0$ and $i \leq \mathbf{I}$.*

PROOF. To obtain the first equality in (3.2), approximate $\chi$ by a sequence of functions $\{f_n\} \subset \mathcal{C}$ such that $0 \leq f_n \uparrow \chi$ and $0 \leq f_n' \uparrow 1_{(0,\infty)}$ as $n \to \infty$, and then use monotone convergence in (3.1), noting property (i) of Definition 3.2. The second equality follows immediately from the definition of $\tau_i$.  □

REMARK. In fact, (3.2) holds also if $\langle \chi, \zeta_i(0) \rangle = \infty$, but then $\langle \chi, \zeta_i(t) \rangle = \infty$ for all $t \geq 0$.

**4. Sequence of systems and fluid limit theorem.** Let $\mathcal{R}$ be a sequence of positive real numbers increasing to infinity. Consider an $\mathcal{R}$-indexed sequence of stochastic models, each defined as in Section 2.3 for the same underlying network structure $(A, C)$ and bandwidth sharing policy $\Lambda$. For each $r \in \mathcal{R}$, there are arrival processes $E_1^r, \ldots, E_\mathbf{I}^r$ with arrival times $\{U_{ik}^r\}_{k=1}^\infty$, $i \leq \mathbf{I}$; there are document sizes $\{v_{1k}^r\}_{k=1}^\infty, \ldots, \{v_{\mathbf{I}k}^r\}_{k=1}^\infty$, with parameters $\vartheta^r$ and $\mu^r$; there is the corresponding measure valued load process $\mathcal{L}^r$; there is an initial condition $\mathcal{Z}^r(0)$; there is a state descriptor $\mathcal{Z}^r$ with auxiliary processes $(Z^r, T^r, U^r, W^r)$ and cumulative service process $S^r(\cdot, \cdot)$. The stochastic elements of each model are defined on a probability space $(\Omega^r, \mathscr{F}^r, \mathbf{P}^r)$ with expectation operator $\mathbf{E}^r$.

4.1. *Scaling.* A fluid scaling (or law of large numbers scaling) is applied to each model in the $\mathcal{R}$-indexed sequence. For each $r \in \mathcal{R}$ and $t \geq s \geq 0$, let

$$\bar{E}^r(t) = \frac{1}{r} E^r(rt), \qquad \bar{S}^r(s,t) = S^r(rs, rt),$$

$$\bar{\mathcal{L}}^r(t) = \frac{1}{r} \mathcal{L}^r(rt), \qquad \bar{\mathcal{L}}^r(s,t) = \frac{1}{r} \mathcal{L}^r(rs, rt),$$

(4.1) $$\bar{\mathcal{Z}}^r(t) = \frac{1}{r} \mathcal{Z}^r(rt), \qquad \bar{Z}^r(t) = \frac{1}{r} Z^r(rt),$$

$$\bar{T}^r(t) = \frac{1}{r} T^r(rt), \qquad \bar{U}^r(t) = \frac{1}{r} U^r(rt),$$

$$\bar{W}^r(t) = \frac{1}{r} W^r(rt).$$

With these definitions, (2.10)–(2.14), and the scaling property of Definition 2.1(iv), we have that for $r \in \mathcal{R}$ and $t \geq 0$,

(4.2) $$\bar{Z}^r(t) = \langle 1, \bar{\mathcal{Z}}^r(t) \rangle,$$

(4.3) $$\bar{T}_i^r(t) = \int_0^t \Lambda_i(\bar{Z}^r(s)) \, ds, \qquad i \leq \mathbf{I},$$

(4.4) $$\bar{U}^r(t) = Ct - A\bar{T}^r(t),$$

(4.5) $$\bar{W}^r(t) = \langle \chi, \bar{\mathcal{Z}}^r(t) \rangle,$$

(4.6) $$\bar{W}^r(t) = \bar{W}^r(0) + \langle \chi, \bar{\mathcal{L}}^r(t) \rangle - \bar{T}^r(t).$$

Also, (2.7) and Definition 2.1(iv) imply that for $r \in \mathcal{R}$ and $[s, t] \subset [0, \infty)$,

(4.7) $$\bar{S}_i^r(s,t) = \int_s^t \frac{\Lambda_i(\bar{Z}^r(u))}{\bar{Z}_i^r(u)} \, du, \qquad i \leq \mathbf{I}.$$

4.2. *Asymptotic assumptions.* In this section, we impose asymptotic assumptions on the $\mathcal{R}$-indexed sequence of models. This is the setting in which our fluid limit result, Theorem 4.1 below, is proved.

Let $\nu = (\nu_1, \ldots, \nu_\mathbf{I})$ be a vector of strictly positive constants and let $\nu(t) = \nu t$ for all $t \geq 0$. Let $\vartheta = (\vartheta_1, \ldots, \vartheta_\mathbf{I})$ be a vector of probability measures in $\mathbf{M}^\mathbf{I}$ satisfying

$$\|\langle 1_{\{0\}}, \vartheta \rangle\| = 0, \tag{4.8}$$

$$\|\langle \chi, \vartheta \rangle\| < \infty. \tag{4.9}$$

For $i \leq \mathbf{I}$, let $\mu_i = \langle \chi, \vartheta_i \rangle^{-1}$ and $\rho_i = \nu_i/\mu_i$. Define $\rho(t) = \rho t$ for all $t \geq 0$. For the sequence of arrival processes, assume that as $r \to \infty$,

$$\bar{E}^r(\cdot) \Rightarrow \nu(\cdot). \tag{4.10}$$

Conditions under which the functional law of large numbers result (4.10) holds are well known. For the sequence of document size distributions, assume that

$$\vartheta^r \xrightarrow{\mathbf{w}} \vartheta \qquad \text{as } r \to \infty, \tag{4.11}$$

$$\{\vartheta_i^r : r \in \mathcal{R}\} \text{ is uniformly integrable for each } i \leq \mathbf{I}. \tag{4.12}$$

Note that (4.11) and (4.12) imply that

$$\mu^r \to \mu \qquad \text{as } r \to \infty. \tag{4.13}$$

For the sequence of fluid scaled initial conditions $\{\bar{\mathcal{Z}}^r(0) : r \in \mathcal{R}\}$, assume that as $r \to \infty$,

$$(\bar{\mathcal{Z}}^r(0), \langle \chi, \bar{\mathcal{Z}}^r(0) \rangle) \Rightarrow (\mathcal{Z}^0, \langle \chi, \mathcal{Z}^0 \rangle), \tag{4.14}$$

where $\mathcal{Z}^0$ is a random vector of measures (taking values in $\mathbf{M}^\mathbf{I}$) satisfying

$$\|\langle \chi, \mathcal{Z}^0 \rangle\| < \infty \qquad \text{a.s.,} \tag{4.15}$$

$$\lim_{\delta \to 0} \mathbf{P}\left(\sup_{x \in \mathbb{R}_+} \|\langle 1_{[x, x+\delta]}, \mathcal{Z}^0 \rangle\| < \varepsilon\right) = 1 \qquad \text{for all } \varepsilon > 0. \tag{4.16}$$

Assumption (4.15) means that the limiting initial workload on each route is finite almost surely; (4.16) is equivalent to the assumption that almost surely, $\mathcal{Z}_i^0$ has no atoms for all $i \leq \mathbf{I}$ (see [7], Lemma A.1).

4.3. *Fluid limit theorem.* The assumptions made so far are now summarized for ease of reference.

(A) *There is a fixed network structure $(A, C)$ and a bandwidth sharing policy $\Lambda$. There is a sequence of stochastic models, each defined as in Section 2.3; there exist a vector of strictly positive constants $\nu$, a vector of probability measures $\vartheta \in \mathbf{M}^\mathbf{I}$, and a random vector of measures $\mathcal{Z}^0$ taking values in $\mathbf{M}^\mathbf{I}$ such that (4.8)–(4.16) hold.*

The following is the main result of the paper.

THEOREM 4.1. *Assume (A). The sequence $\{(\bar{\mathcal{Z}}^r, \bar{Z}^r, \bar{T}^r, \bar{U}^r, \bar{W}^r)\}$ is **C**-tight, and each weak limit point $(\mathcal{Z}, Z, T, U, W)$ is such that almost surely, $\mathcal{Z}$ is a fluid model solution with auxiliary functions $(Z, T, U, W)$ for the data $(A, C, \Lambda, \nu, \vartheta)$, where $W(t)$ is finite for all $t \geq 0$.*

**5. Proof of Theorem 4.1.** The proof has several stages. Section 5.1 contains a functional law of large numbers result for the measure valued load processes $\{\mathcal{L}^r\}$. This result follows from the assumptions imposed on the stochastic primitives. Section 5.2 derives two dynamic equations satisfied by the fluid scaled state descriptors $\{\bar{\mathcal{Z}}^r\}$, as well as several related bounds. Section 5.3 establishes a compact containment property, and Sections 5.4 and 5.5 establish control of oscillations for the state descriptors. These properties are combined in Section 5.6 to prove the tightness claim of Theorem 4.1, and properties of weak limit points are derived in Section 5.7. We assume (A) throughout this entire section.

The general strategy outlined above is similar to that in [7]. However, the model studied here presents the additional complication of multiple routes that interact with each other via the bandwidth sharing policy $\Lambda$. In particular, the numerator in the first integral term of (3.1) is a function of the current state of the whole system, as opposed to a constant as is the case in the analogous equation in [7]. This requires additional care to carry out the analysis. A key difference in the present proof is in verifying (at various stages along the way), that the assumptions imposed by Definition 2.1 on the more general function $\Lambda$ are sufficient to allow the above strategy to go through. Furthermore, [7] focused only on a heavily loaded single server queue and its critical fluid limit. Here, we have a network of resources and there is no a priori assumption on the system load, that is, the traffic intensity parameters $\rho_i$ are unrestricted in $(0, \infty)$. This results in a more subtle fluid model and limit proof (see Section 5.7) related to the treatment of times when fluid queue lengths become zero.

5.1. *Limit of the primitive load processes.* Recall that $\nu(t) = \nu t$, and $\rho(t) = \rho t$ for all $t \geq 0$.

THEOREM 5.1. *As $r \to \infty$,*

$$(5.1) \qquad (\bar{\mathcal{L}}^r(\cdot), \langle \chi, \bar{\mathcal{L}}^r(\cdot) \rangle) \Rightarrow (\nu(\cdot)\vartheta, \rho(\cdot)).$$

The proof of this theorem is a straightforward application of a functional law of large numbers. For completeness, a proof is given in the Appendix.

5.2. *Dynamic equations.* Fix $r \in \mathcal{R}$. For each route $i \leq \mathbf{I}$, a dynamic equation satisfied by the component $\bar{\mathcal{Z}}_i^r(\cdot)$ of the fluid scaled state descriptor is the starting point for much of our subsequent analysis. The equation results, after some simplification, from substituting the definition of the residual document sizes (2.8) into (2.9). Almost surely, for all Borel measurable $f : \mathbb{R}_+ \to \mathbb{R}$, all $i \leq \mathbf{I}$, and all $t \geq s \geq 0$,

$$\langle f, \mathcal{Z}_i^r(t) \rangle = \langle f(\cdot - S_i^r(s,t)), \mathcal{Z}_i^r(s) \rangle + \sum_{k=E_i^r(s)+1}^{E_i^r(t)} f(v_{ik}^r - S_i^r(U_{ik}^r, t)).$$

Recall that $f$ is always extended to be zero on $(-\infty, 0)$ so that $f(\cdot - x)$ is well defined on $\mathbb{R}_+$ for all $x \geq 0$. Applying the fluid scaling (4.1) produces

$$(5.2) \quad \langle f, \bar{\mathcal{Z}}_i^r(t) \rangle = \langle f(\cdot - \bar{S}_i^r(s,t)), \bar{\mathcal{Z}}_i^r(s) \rangle + \frac{1}{r} \sum_{k=r\bar{E}_i^r(s)+1}^{r\bar{E}_i^r(t)} f(v_{ik}^r - \bar{S}_i^r(U_{ik}^r r^{-1}, t)).$$

This equation yields several estimates that will be used frequently. If $f$ is nonnegative and nondecreasing, then using the bound $\sup_{x \in \mathbb{R}_+} f(\cdot - x) \leq f(\cdot)$ in (5.2) yields

$$(5.3) \quad \begin{aligned} \langle f, \bar{\mathcal{Z}}_i^r(t) \rangle &\leq \langle f(\cdot - \bar{S}_i^r(s,t)), \bar{\mathcal{Z}}_i^r(s) \rangle + \langle f, \bar{\mathcal{L}}_i^r(s,t) \rangle \\ &\leq \langle f, \bar{\mathcal{Z}}_i^r(s) \rangle + \langle f, \bar{\mathcal{L}}_i^r(s,t) \rangle. \end{aligned}$$

If $f$ is bounded, then (5.2) implies that

$$(5.4) \quad \begin{aligned} \langle f, \bar{\mathcal{Z}}_i^r(t) \rangle &\leq \langle f(\cdot - \bar{S}_i^r(s,t)), \bar{\mathcal{Z}}_i^r(s) \rangle + \|f\|_\infty \langle 1, \bar{\mathcal{L}}_i^r(s,t) \rangle \\ &\leq \|f\|_\infty \langle 1, \bar{\mathcal{Z}}_i^r(s) \rangle + \|f\|_\infty \langle 1, \bar{\mathcal{L}}_i^r(s,t) \rangle. \end{aligned}$$

By ignoring the sum in (5.2), we obtain for any nonnegative $f$ that

$$(5.5) \quad \langle f(\cdot - \bar{S}_i^r(s,t)), \bar{\mathcal{Z}}_i^r(s) \rangle \leq \langle f, \bar{\mathcal{Z}}_i^r(t) \rangle.$$

An alternative dynamic equation to (5.2), that is satisfied by $\bar{\mathcal{Z}}_i^r(\cdot)$ on certain time intervals, will be used when passing to the limit as $r \to \infty$. This equation is a prelimit analogue of the (3.1) satisfied by fluid model solutions. It is derived from (5.2) and is written in terms of projections against functions $f$ in the more restrictive class

$$(5.6) \quad \mathcal{C}_c = \{ f \in \mathcal{C} : f \text{ has compact support in } \mathbb{R}_+ \}.$$

Note that for $f \in \mathcal{C}_c$, the derivative $f'$ has compact support and $\|f'\|_\infty < \infty$. The proof of the following result appears in the Appendix.

LEMMA 5.2. *Fix $r \in \mathcal{R}$. Almost surely, for all $i \leq \mathbf{I}$, all $f \in \mathcal{C}_c$, and all finite time intervals $[s,t] \subset [0,\infty)$ satisfying $\inf_{u \in [s,t]} \bar{Z}_i^r(u) > 0$, we have*

$$\langle f, \bar{\mathcal{Z}}_i^r(t) \rangle = \langle f, \bar{\mathcal{Z}}_i^r(s) \rangle - \int_s^t \langle f', \bar{\mathcal{Z}}_i^r(u) \rangle \frac{\Lambda_i(\bar{Z}^r(u))}{\bar{Z}_i^r(u)} \, du$$
(5.7)
$$+ \langle f, \bar{\mathcal{L}}_i^r(t) \rangle - \langle f, \bar{\mathcal{L}}_i^r(s) \rangle.$$

5.3. *Compact containment.* In this section, we establish the first of the two main conditions used in proving tightness.

LEMMA 5.3. *Let $T > 0$ and $\eta > 0$. There exists a compact set $\mathbf{K} \subset \mathbf{M}^{\mathbf{I}}$ such that*

$$\liminf_{r \to \infty} \mathbf{P}^r(\bar{\mathcal{Z}}^r(t) \in \mathbf{K} \text{ for all } t \in [0,T]) \geq 1 - \eta.$$

PROOF. By (4.15) and since $\|\langle 1, \mathcal{Z}^0 \rangle\| < \infty$ almost surely, there exists an $M > 0$ such that

(5.8) $$\mathbf{P}(\|\langle 1, \mathcal{Z}^0 \rangle\| \vee \|\langle \chi, \mathcal{Z}^0 \rangle\| \geq M) \leq \eta.$$

Since $\xi \mapsto \langle 1, \xi \rangle$ is a continuous $\mathbb{R}_+^{\mathbf{I}}$-valued function on $\mathbf{M}^{\mathbf{I}}$, assumption (4.14) and the continuous mapping theorem imply that

(5.9) $$(\langle 1, \bar{\mathcal{Z}}^r(0) \rangle, \langle \chi, \bar{\mathcal{Z}}^r(0) \rangle) \Rightarrow (\langle 1, \mathcal{Z}^0 \rangle, \langle \chi, \mathcal{Z}^0 \rangle) \quad \text{as } r \to \infty.$$

The set $\{(z,w) \in \mathbb{R}_+^{\mathbf{I}} \times \mathbb{R}_+^{\mathbf{I}} : \|z\| \vee \|w\| < M\}$ is open, so by (5.8), (5.9) and the Portmanteau theorem,

$$\liminf_{r \to \infty} \mathbf{P}^r(\|\langle 1, \bar{\mathcal{Z}}^r(0) \rangle\| \vee \|\langle \chi, \bar{\mathcal{Z}}^r(0) \rangle\| < M)$$
(5.10)
$$\geq \mathbf{P}(\|\langle 1, \mathcal{Z}^0 \rangle\| \vee \|\langle \chi, \mathcal{Z}^0 \rangle\| < M) \geq 1 - \eta.$$

For each $r \in \mathcal{R}$, let $\Omega_1^r$ be the event in the left-hand side of (5.10) and define

$$\Omega_2^r = \{\|\langle 1, \bar{\mathcal{L}}^r(T) \rangle\| \vee \|\langle \chi, \bar{\mathcal{L}}^r(T) \rangle\| < K\},$$

where $K = (\|\nu T\| \vee \|\rho T\|) + 1$. By Theorem 5.1,

$$(\langle 1, \bar{\mathcal{L}}^r(T) \rangle, \langle \chi, \bar{\mathcal{L}}^r(T) \rangle) \Rightarrow (\nu T, \rho T) \quad \text{as } r \to \infty.$$

So $\liminf_{r \to \infty} \mathbf{P}^r(\Omega_2^r) = 1$ by the choice of $K$. For each $r \in \mathcal{R}$, let $\Omega_3^r$ be a full probability event on which the dynamic equation (5.2) holds. Then

(5.11) $$\liminf_{r \to \infty} \mathbf{P}^r(\Omega_1^r \cap \Omega_2^r \cap \Omega_3^r) \geq 1 - \eta.$$

Let $\mathbf{K}$ be the closure in $\mathbf{M}^{\mathbf{I}}$ of the set $\{\xi \in \mathbf{M}^{\mathbf{I}} : \|\langle 1, \xi \rangle\| \vee \|\langle \chi, \xi \rangle\| \leq M + K\}$. The set $\mathbf{K}$ is compact by [11], Theorem 15.7.5. Fix $r \in \mathcal{R}$ and an outcome $\omega \in \Omega_1^r \cap \Omega_2^r \cap \Omega_3^r$; assume for the rest of the proof that all random objects are

evaluated at this $\omega$. Fix $t \in [0, T]$; by (5.11), it suffices to show that $\bar{\mathcal{Z}}^r(t) \in \mathbf{K}$. The dynamic equation bound (5.4) and the definition of $\Omega^r_1 \cap \Omega^r_2 \cap \Omega^r_3$ imply that

$$\max_{i \leq \mathbf{I}} \langle 1, \bar{\mathcal{Z}}^r_i(t) \rangle \leq \max_{i \leq \mathbf{I}} \{ \langle 1, \bar{\mathcal{Z}}^r_i(0) \rangle + \langle 1, \bar{\mathcal{L}}^r_i(t) \rangle \}$$

(5.12)
$$\leq \max_{i \leq \mathbf{I}} \{ \langle 1, \bar{\mathcal{Z}}^r_i(0) \rangle + \langle 1, \bar{\mathcal{L}}^r_i(T) \rangle \}$$

$$\leq M + K.$$

Similarly, the dynamic equation bound (5.3) implies that

$$\max_{i \leq \mathbf{I}} \langle \chi, \bar{\mathcal{Z}}^r_i(t) \rangle \leq \max_{i \leq \mathbf{I}} \{ \langle \chi, \bar{\mathcal{Z}}^r_i(0) \rangle + \langle \chi, \bar{\mathcal{L}}^r_i(t) \rangle \}$$

(5.13)
$$\leq \max_{i \leq \mathbf{I}} \{ \langle \chi, \bar{\mathcal{Z}}^r_i(0) \rangle + \langle \chi, \bar{\mathcal{L}}^r_i(T) \rangle \}$$

$$\leq M + K.$$

Combining (5.12) and (5.13) with (5.11) completes the proof. □

5.4. *Asymptotic regularity near zero.* Over any finite time interval, with arbitrarily high probability as $r \to \infty$, the fluid scaled state descriptor $\bar{\mathcal{Z}}^r_i(\cdot)$ for route $i$ puts arbitrarily small mass on a sufficiently small neighborhood of zero. This is proved in the following lemma, and is a key ingredient for establishing an oscillation property in the next section.

LEMMA 5.4. *Let $T > 0$. For each $\varepsilon, \eta \in (0, 1)$, there exists an $a > 0$ such that*

(5.14)
$$\liminf_{r \to \infty} \mathbf{P}^r \left( \sup_{t \in [0,T]} \| \langle 1_{[0,a]}, \bar{\mathcal{Z}}^r(t) \rangle \| \leq \varepsilon \right) \geq 1 - \eta.$$

PROOF. Fix $\varepsilon, \eta \in (0, 1)$. The proof consists of several steps. The first three steps are concerned with defining four high probability events $\Omega^r_1, \Omega^r_2, \Omega^r_3, \Omega^r_4$. Steps four and five supply the desired bound (in two parts) on the intersection of these events.

*Step* 1. By (4.16), there exists $b > 0$ such that

(5.15)
$$\mathbf{P} \left( \sup_{n \in \mathbb{N}} \| \langle 1_{[(n-1)b, nb]}, \mathcal{Z}^0 \rangle \| < \frac{\varepsilon}{4} \right) \geq 1 - \frac{\eta}{2}.$$

Let $\mathbf{B} = \{ \xi \in \mathbf{M}^\mathbf{I} : \sup_{n \in \mathbb{N}} \| \langle 1_{[(n-1)b, nb]}, \xi \rangle \| < \varepsilon/4 \}$ and suppose that $\xi \in \mathbf{B}$ and $\{ \xi^k \} \subset \mathbf{M}^\mathbf{I}$ satisfy $\xi^k \xrightarrow{\mathbf{w}} \xi$. Choose $L \in \mathbb{N}$ large enough so that $\| \langle 1_{[Lb, \infty)}, \xi \rangle \| < \varepsilon/4$. Since $\xi^k_i \xrightarrow{\mathbf{w}} \xi_i$ for each $i \leq \mathbf{I}$, the Portmanteau theorem implies that

$$\limsup_{k \to \infty} \sup_{n \in \mathbb{N}} \| \langle 1_{[(n-1)b, nb]}, \xi^k \rangle \|$$

$$\leq \limsup_{k\to\infty} \left( \max_{n\leq L} \|\langle 1_{[(n-1)b,nb]}, \xi^k\rangle\| \vee \|\langle 1_{[Lb,\infty)}, \xi^k\rangle\| \right)$$

$$\leq \max_{n\leq L} \|\langle 1_{[(n-1)b,nb]}, \xi\rangle\| \vee \|\langle 1_{[Lb,\infty)}, \xi\rangle\| < \frac{\varepsilon}{4}.$$

So $\xi^k \in \mathbf{B}$ for sufficiently large $k$, which implies that $\mathbf{B} \subset \mathbf{M}^{\mathbf{I}}$ is open. We deduce from (4.14) and the Portmanteau theorem that

$$\liminf_{r\to\infty} \mathbf{P}^r \left( \sup_{x\in\mathbb{R}_+} \|\langle 1_{[x,x+b]}, \bar{\mathcal{Z}}^r(0)\rangle\| < \frac{\varepsilon}{2} \right)$$

(5.16)
$$\geq \liminf_{r\to\infty} \mathbf{P}^r \left( \sup_{n\in\mathbb{N}} \|\langle 1_{[(n-1)b,nb]}, \bar{\mathcal{Z}}^r(0)\rangle\| < \frac{\varepsilon}{4} \right)$$

$$\geq \mathbf{P} \left( \sup_{n\in\mathbb{N}} \|\langle 1_{[(n-1)b,nb]}, \mathcal{Z}^0\rangle\| < \frac{\varepsilon}{4} \right).$$

Combining (5.16) with (5.15) yields

(5.17) $$\liminf_{r\to\infty} \mathbf{P}^r \left( \sup_{x\in\mathbb{R}_+} \|\langle 1_{[x,x+b]}, \bar{\mathcal{Z}}^r(0)\rangle\| < \frac{\varepsilon}{2} \right) \geq 1 - \frac{\eta}{2}.$$

Let $\Omega_1^r$ be the event in the left-hand side of (5.17).

*Step* 2. By Lemma 5.3, there exists a compact set $\mathbf{K} \subset \mathbf{M}^{\mathbf{I}}$ such that

(5.18) $$\liminf_{r\to\infty} \mathbf{P}^r(\bar{\mathcal{Z}}^r(t) \in \mathbf{K} \text{ for all } t \in [0,T]) \geq 1 - \frac{\eta}{2}.$$

Since $\mathbf{K}$ is compact, there exists $M < \infty$ such that

(5.19) $$\sup_{\xi\in\mathbf{K}} \|\langle 1, \xi\rangle\| \leq M.$$

Let $\Omega_2^r$ be the event in the left-hand side of (5.18).

*Step* 3. By Lemma 2.2, there exists $c > 0$ such that for each $i \leq \mathbf{I}$,

(5.20) $$\Lambda_i(z) \geq c \qquad \text{on } \{z \in \mathbb{R}_+^{\mathbf{I}} : z_i \geq \varepsilon/8, \|z\| \leq M\}.$$

Let $\delta = \varepsilon(12\|\nu\|)^{-1} \wedge T$ and let $a = \delta c(2M)^{-1} \wedge b$. Choose $N \in \mathbb{N}$ large enough so that

(5.21) $$Na > a + T\|C\|\frac{8}{\varepsilon}.$$

Let $I_0 = \varnothing$ and, for each $n \in \mathbb{N}$, define $I_n = [(n-1)a, na)$ and choose $g_n \in \mathbf{C}_b(\mathbb{R}_+)$ satisfying $1_{I_n} \leq g_n \leq 1_{I_{n-1}\cup I_n\cup I_{n+1}}$. Then since $\vartheta$ is a vector of probability measures,

(5.22) $$\max_{i\leq\mathbf{I}} \sum_{n=1}^{\infty} \langle g_n, \vartheta_i\rangle \leq \max_{i\leq\mathbf{I}} \sum_{n=1}^{\infty} \langle 1_{I_{n-1}\cup I_n\cup I_{n+1}}, \vartheta_i\rangle \leq 3.$$

For notational convenience, let $g_0 \equiv 1$. For each $n \in \mathbb{N} \cup \{0\}$, the $\mathbb{R}_+^{\mathbf{I}}$-valued map $\xi \mapsto \langle g_n, \xi \rangle$ is continuous on $\mathbf{M}^{\mathbf{I}}$. So for each such $n$, Theorem 5.1 and the continuous mapping theorem yield

$$(5.23) \qquad \langle g_n, \bar{\mathcal{L}}^r(\cdot) \rangle \Rightarrow \nu(\cdot) \langle g_n, \vartheta \rangle \qquad \text{as } r \to \infty.$$

The limit in (5.23) is a deterministic and continuous function taking values in $\mathbb{R}_+^{\mathbf{I}}$. So, the convergence is uniform on compact time intervals in probability, and occurs jointly for all $n = 0, \ldots, N$. Therefore,

$$(5.24) \quad \liminf_{r \to \infty} \mathbf{P}^r \left( \max_{n=0,\ldots,N} \sup_{t \in [0,T]} \|\langle g_n, \bar{\mathcal{L}}^r(t) \rangle - \nu t \langle g_n, \vartheta \rangle \| \leq \frac{\varepsilon}{8N} \right) = 1.$$

Let $\Omega_3^r$ be the event in (5.24) and let $\Omega_4^r$ be a full probability event on which (5.2) holds. Define $\Omega_0^r = \Omega_1^r \cap \Omega_2^r \cap \Omega_3^r \cap \Omega_4^r$. By (5.17), (5.18) and (5.24),

$$\liminf_{r \to \infty} \mathbf{P}^r (\Omega_0^r) \geq 1 - \eta.$$

Thus, setting $\Omega_*^r$ equal to the event in (5.14), it suffices to show that $\Omega_0^r \subset \Omega_*^r$ for each $r \in \mathcal{R}$. To this end, fix $r \in \mathcal{R}$, $\omega \in \Omega_0^r$, $t \in [0,T]$, and $i \leq \mathbf{I}$; assume for the rest of the proof that all random objects are evaluated at this $\omega$. It suffices to show that

$$(5.25) \qquad \langle 1_{[0,a]}, \bar{\mathcal{Z}}_i^r(t) \rangle \leq \varepsilon.$$

*Step* 4. Define the random time

$$(5.26) \qquad \tau = \sup \left\{ s \leq t : \bar{Z}_i^r(s) \leq \frac{\varepsilon}{8} \right\},$$

where $\sup \varnothing = 0$. We first show that

$$(5.27) \qquad \sup_{x \in \mathbb{R}_+} \langle 1_{[0,a]}(\cdot - x), \bar{\mathcal{Z}}_i^r(\tau) \rangle \leq \frac{\varepsilon}{2}.$$

If $\tau = 0$, this follows from the definition of $\Omega_1^r$, since $a \leq b$. If $\tau > 0$, then the definition of $\tau$ implies the existence of $\tilde{\tau} \in [(\tau - \delta)^+, \tau]$ such that $\langle 1, \bar{\mathcal{Z}}_i^r(\tilde{\tau}) \rangle = \bar{Z}_i^r(\tilde{\tau}) \leq \varepsilon/8$. By the dynamic equation bound (5.4) and the definition of $\tilde{\tau}$,

$$\sup_{x \in \mathbb{R}_+} \langle 1_{[0,a]}(\cdot - x), \bar{\mathcal{Z}}_i^r(\tau) \rangle \leq \langle 1, \bar{\mathcal{Z}}_i^r(\tilde{\tau}) \rangle + \langle 1, \bar{\mathcal{L}}_i^r(\tilde{\tau}, \tau) \rangle$$

$$\leq \frac{\varepsilon}{8} + \langle 1, \bar{\mathcal{L}}_i^r(\tau) \rangle - \langle 1, \bar{\mathcal{L}}_i^r(\tilde{\tau}) \rangle.$$

Applying the definition of $\Omega_3^r$ and noting that $g_0 \equiv 1$, we obtain

$$\sup_{x \in \mathbb{R}_+} \langle 1_{[0,a]}(\cdot - x), \bar{\mathcal{Z}}_i^r(\tau) \rangle \leq \frac{\varepsilon}{8} + \nu_i (\tau - \tilde{\tau}) \langle 1, \vartheta_i \rangle + \frac{\varepsilon}{4N} \leq \frac{\varepsilon}{8} + \delta \|\nu\| + \frac{\varepsilon}{4},$$

which implies (5.27) by the choice of $\delta$.

*Step* 5. Note that if $\tau = t$, then (5.25) follows directly from (5.27); so assume that $t > \tau$. For all $s \in (\tau, t]$, $\bar{Z}_i^r(s) > \varepsilon/8$ and $\|\bar{Z}^r(s)\| \leq M$ by (5.19) and the definition of $\Omega_2^r$. So, (5.20) implies that

$$\inf_{s \in (\tau, t]} \Lambda_i(\bar{Z}^r(s)) \geq c. \tag{5.28}$$

Using (5.2) and (5.27),

$$\langle 1_{[0,a]}, \bar{\mathcal{Z}}_i^r(t) \rangle = \langle 1_{[0,a]}(\cdot - \bar{S}_i^r(\tau, t)), \bar{\mathcal{Z}}_i^r(\tau) \rangle$$

$$+ \frac{1}{r} \sum_{k=r\bar{E}_i^r(\tau)+1}^{r\bar{E}_i^r(t)} 1_{[0,a]}(v_{ik}^r - \bar{S}_i^r(U_{ik}^r r^{-1}, t)) \tag{5.29}$$

$$\leq \frac{\varepsilon}{2} + \sum_{n=1}^{\infty} \frac{1}{r} \sum_{k=r\bar{E}_i^r(\tau)+1}^{r\bar{E}_i^r(t)} 1_{I_n}(v_{ik}^r) 1_{[0,a]}(v_{ik}^r - \bar{S}_i^r(U_{ik}^r r^{-1}, t)).$$

Consider a flow $k$ such that $U_{ik}^r r^{-1} \in (\tau, t]$ and $v_{ik}^r \in I_n$ for $n > N$. Then $v_{ik}^r \geq Na > a + T\|C\|8\varepsilon^{-1}$. Since $\bar{Z}_i^r(s) > \varepsilon/8$ for $s \in (\tau, t]$,

$$\bar{S}_i^r(U_{ik}^r r^{-1}, t) = \int_{U_{ik}^r r^{-1}}^{t} \frac{\Lambda_i(\bar{Z}^r(s))}{\bar{Z}_i^r(s)} ds \leq T\|C\| \frac{8}{\varepsilon}.$$

Thus, $(v_{ik}^r - \bar{S}_i^r(U_{ik}^r r^{-1}, t)) > a$ and so $1_{[0,a]}(v_{ik}^r - \bar{S}_i^r(U_{ik}^r r^{-1}, t)) = 0$. We deduce from (5.29) that

$$\langle 1_{[0,a]}, \bar{\mathcal{Z}}_i^r(t) \rangle \leq \frac{\varepsilon}{2} + \sum_{n=1}^{N} \frac{1}{r} \sum_{k=r\bar{E}_i^r(\tau)+1}^{r\bar{E}_i^r(t)} 1_{I_n}(v_{ik}^r) 1_{[0,a]}(v_{ik}^r - \bar{S}_i^r(U_{ik}^r r^{-1}, t)). \tag{5.30}$$

Consider two flows $k < l$ satisfying $U_{ik}^r r^{-1}, U_{il}^r r^{-1} \in (\tau, t]$ and $v_{ik}^r, v_{il}^r \in I_n$ for some $n = 1, \ldots, N$. If $U_{il}^r r^{-1} - U_{ik}^r r^{-1} \geq \delta$, then by the definition of $\Omega_2^r$, (5.19), (5.28) and the definition of $a$,

$$\bar{S}_i^r(U_{ik}^r r^{-1}, t) - \bar{S}_i^r(U_{il}^r r^{-1}, t) = \int_{U_{ik}^r r^{-1}}^{U_{il}^r r^{-1}} \frac{\Lambda_i(\bar{Z}^r(s))}{\bar{Z}_i^r(s)} ds$$

$$\geq (U_{il}^r r^{-1} - U_{ik}^r r^{-1}) \frac{c}{M}$$

$$\geq \frac{\delta c}{M} \geq 2a.$$

Consequently,

$$(v_{il}^r - \bar{S}_i^r(U_{il}^r r^{-1}, t)) - (v_{ik}^r - \bar{S}_i^r(U_{ik}^r r^{-1}, t)) \geq 2a + v_{il}^r - v_{ik}^r > 2a - a = a$$

and so at most one of

$$1_{[0,a]}(v_{ik}^r - \bar{S}_i^r(U_{ik}^r r^{-1}, t)) \quad \text{and} \quad 1_{[0,a]}(v_{il}^r - \bar{S}_i^r(U_{il}^r r^{-1}, t))$$

is nonzero. This implies that flows arriving to route $i$ during $(\tau, t]$ with document sizes in $I_n$ and residual document sizes at time $t$ in $[0, a]$ must all arrive during some time interval of length less than $\delta$. That is, for each $n = 1, \ldots, N$, there exists an interval $(s_n, s_n + \delta_n] \subset (\tau, t]$, with $\delta_n < \delta$, such that $U_{ik}^r r^{-1} \in (\tau, t]$, $v_{ik}^r \in I_n$, and $v_{ik}^r - \bar{S}_i^r(U_{ik}^r r^{-1}, t) \in [0, a]$, implies $U_{ik}^r r^{-1} \in (s_n, s_n + \delta_n]$. Combining this fact with (5.30) yields

$$\langle 1_{[0,a]}, \bar{\mathcal{Z}}_i^r(t) \rangle \leq \frac{\varepsilon}{2} + \sum_{n=1}^{N} \sup_{s \in [0, T-\delta]} \frac{1}{r} \sum_{k = r\bar{E}_i^r(s)+1}^{r\bar{E}_i^r(s+\delta)} 1_{I_n}(v_{ik}^r).$$

Bound $1_{I_n}$ by $g_n$ and rewrite the above to obtain

$$\langle 1_{[0,a]}, \bar{\mathcal{Z}}_i^r(t) \rangle \leq \frac{\varepsilon}{2} + \sum_{n=1}^{N} \sup_{s \in [0, T-\delta]} (\langle g_n, \bar{\mathcal{L}}_i^r(s+\delta) \rangle - \langle g_n, \bar{\mathcal{L}}_i^r(s) \rangle).$$

Applying the definition of $\Omega_3^r$ and (5.22), we obtain

$$\langle 1_{[0,a]}, \bar{\mathcal{Z}}_i^r(t) \rangle \leq \frac{\varepsilon}{2} + \sum_{n=1}^{N} \left( \nu_i \delta \langle g_n, \vartheta_i \rangle + \frac{\varepsilon}{4N} \right) \leq \frac{3\varepsilon}{4} + \delta 3 \|\nu\|.$$

By the choice of $\delta$, the right-hand side is bounded above by $\varepsilon$.  □

5.5. *Oscillations.* This section contains an oscillation bound used in proving tightness.

DEFINITION 5.5.  Let $T > 0$ and $\delta \in [0, T]$. For each $\zeta(\cdot) \in \mathbf{D}([0, \infty), \mathbf{M}^{\mathbf{I}})$, define a modulus of continuity on $[0, T]$ by

(5.31) $$\mathbf{w}_T(\zeta(\cdot), \delta) = \sup_{\substack{s,t \in [0,T] \\ |t-s| < \delta}} \mathbf{d}_{\mathbf{I}}[\zeta(s), \zeta(t)].$$

LEMMA 5.6.  *Let $T > 0$. For all $\varepsilon, \eta \in (0, 1)$, there exists $\delta > 0$ such that*

(5.32) $$\liminf_{r \to \infty} \mathbf{P}^r(\mathbf{w}_T(\bar{\mathcal{Z}}^r(\cdot), \delta) \leq \varepsilon) \geq 1 - \eta.$$

PROOF.  Fix $\varepsilon, \eta \in (0, 1)$. By Lemma 5.4, with $\varepsilon/2$ in place of $\varepsilon$, there exists $a > 0$ such that

(5.33) $$\liminf_{r \to \infty} \mathbf{P}^r \left( \sup_{t \in [0,T]} \|\langle 1_{[0,a]}, \bar{\mathcal{Z}}^r(t) \rangle\| \leq \frac{\varepsilon}{2} \right) \geq 1 - \eta.$$

Let $\Omega_1^r$ be the event in (5.33) and let $\delta = \min\{\varepsilon(\varepsilon \wedge a)(4\|C\|)^{-1}, \varepsilon(4\|\nu\|)^{-1}\}$. Since $\xi \mapsto \langle 1, \xi \rangle$ is continuous on $\mathbf{M}^{\mathbf{I}}$, Theorem 5.1 implies that $\langle 1, \bar{\mathcal{L}}^r(\cdot) \rangle \Rightarrow \nu(\cdot)$ as $r \to \infty$. Thus,

(5.34) $$\liminf_{r \to \infty} \mathbf{P}^r \left( \sup_{t \in [0,T]} \|\langle 1, \bar{\mathcal{L}}^r(t) \rangle - \nu t\| \leq \frac{\varepsilon}{8} \right) = 1.$$

Let $\Omega_2^r$ be the event in (5.34) and let $\Omega_3^r$ be a full probability event on which (5.2) holds. By (5.33) and (5.34),

(5.35) $$\liminf_{r\to\infty} \mathbf{P}^r(\Omega_1^r \cap \Omega_2^r \cap \Omega_3^r) \geq 1 - \eta.$$

Fix $r \in \mathcal{R}$ and an outcome $\omega \in \Omega_1^r \cap \Omega_2^r \cap \Omega_3^r$; assume for the rest of the proof that all random objects are evaluated at this $\omega$. Fix $i \leq \mathbf{I}$ and $s, t \in [0, T]$ with $s \leq t$ and $t - s < \delta$. By (5.35), Definition 5.5 and (1.2) it suffices to show that

$$\mathbf{d}[\bar{\mathcal{Z}}_i^r(s), \bar{\mathcal{Z}}_i^r(t)] \leq \varepsilon.$$

Let $B \subset \mathbb{R}_+$ be closed. By (1.1), it suffices to show the two inequalities,

(5.36) $$\langle 1_B, \bar{\mathcal{Z}}_i^r(s) \rangle \leq \langle 1_{B^\varepsilon}, \bar{\mathcal{Z}}_i^r(t) \rangle + \varepsilon,$$

(5.37) $$\langle 1_B, \bar{\mathcal{Z}}_i^r(t) \rangle \leq \langle 1_{B^\varepsilon}, \bar{\mathcal{Z}}_i^r(s) \rangle + \varepsilon.$$

To show (5.36), use the definition of $\Omega_1^r$ to write

(5.38) $$\begin{aligned}\langle 1_B, \bar{\mathcal{Z}}_i^r(s) \rangle &\leq \langle 1_{[0,a]}, \bar{\mathcal{Z}}_i^r(s) \rangle + \langle 1_{B \cap (a,\infty)}, \bar{\mathcal{Z}}_i^r(s) \rangle \\ &\leq \frac{\varepsilon}{2} + \langle 1_{B \cap (a,\infty)}, \bar{\mathcal{Z}}_i^r(s) \rangle.\end{aligned}$$

Let $I = \{u \in [s,t] : \bar{Z}_i^r(u) < \varepsilon/2\}$. Suppose that $I = \varnothing$. Then $\bar{Z}_i^r(u) \geq \varepsilon/2$ for all $u \in [s,t]$. So, by (2.1) and the definition of $\delta$,

(5.39) $$\bar{S}_i^r(s,t) = \int_s^t \frac{\Lambda_i(\bar{Z}^r(u))}{\bar{Z}_i^r(u)}\, du \leq \delta \frac{2\|C\|}{\varepsilon} \leq \frac{\varepsilon \wedge a}{2} < \varepsilon \wedge a.$$

Consequently, $x \in B \cap (a, \infty)$ implies $x - \bar{S}_i^r(s,t) \in B^\varepsilon$, and so

$$B \cap (a, \infty) \subset B^\varepsilon + \bar{S}_i^r(s,t).$$

We deduce from (5.38) that

$$\langle 1_B, \bar{\mathcal{Z}}_i^r(s) \rangle \leq \frac{\varepsilon}{2} + \langle 1_{B^\varepsilon + \bar{S}_i^r(s,t)}, \bar{\mathcal{Z}}_i^r(s) \rangle = \frac{\varepsilon}{2} + \langle 1_{B^\varepsilon}(\cdot - \bar{S}_i^r(s,t)), \bar{\mathcal{Z}}_i^r(s) \rangle.$$

Apply the dynamic equation bound (5.5) to obtain

(5.40) $$\langle 1_B, \bar{\mathcal{Z}}_i^r(s) \rangle \leq \frac{\varepsilon}{2} + \langle 1_{B^\varepsilon}, \bar{\mathcal{Z}}_i^r(t) \rangle.$$

Now, suppose $I \neq \varnothing$ and let $\tau = \inf I$. Then by right continuity of $\bar{Z}_i^r(\cdot)$,

(5.41) $$\bar{Z}_i^r(\tau) \leq \frac{\varepsilon}{2}.$$

Since $\bar{Z}_i^r(u) \geq \varepsilon/2$ for all $u \in [s, \tau)$,

(5.42) $$\bar{S}_i^r(s,\tau) = \int_s^\tau \frac{\Lambda_i(\bar{Z}^r(u))}{\bar{Z}_i^r(u)}\, du \leq \delta \frac{2\|C\|}{\varepsilon} \leq a.$$

By (5.38) and (5.42),

$$\langle 1_B, \bar{\mathcal{Z}}_i^r(s)\rangle \leq \frac{\varepsilon}{2} + \langle 1_{(a,\infty)}, \bar{\mathcal{Z}}_i^r(s)\rangle$$

$$\leq \frac{\varepsilon}{2} + \langle 1_{[\bar{S}_i^r(s,\tau),\infty)}, \bar{\mathcal{Z}}_i^r(s)\rangle$$

$$= \frac{\varepsilon}{2} + \langle 1_{[0,\infty)}(\cdot - \bar{S}_i^r(s,\tau)), \bar{\mathcal{Z}}_i^r(s)\rangle.$$

Apply the dynamic equation bound (5.5) to obtain

(5.43) $$\langle 1_B, \bar{\mathcal{Z}}_i^r(s)\rangle \leq \frac{\varepsilon}{2} + \langle 1, \bar{\mathcal{Z}}_i^r(\tau)\rangle \leq \varepsilon.$$

So (5.36) follows because either (5.40) or (5.43) holds.

To show (5.37), note that by definition of $\Omega_2^r$ and $\delta$,

(5.44) $$\langle 1, \bar{\mathcal{L}}_i^r(s,t)\rangle = \langle 1, \bar{\mathcal{L}}_i^r(t)\rangle - \langle 1, \bar{\mathcal{L}}_i^r(s)\rangle \leq \nu_i(t-s) + \frac{\varepsilon}{4} \leq \|\nu\|\delta + \frac{\varepsilon}{4} \leq \frac{\varepsilon}{2}.$$

By the first inequality in (5.4) and (5.44),

(5.45) $$\langle 1_B, \bar{\mathcal{Z}}_i^r(t)\rangle \leq \langle 1_B(\cdot - \bar{S}_i^r(s,t)), \bar{\mathcal{Z}}_i^r(s)\rangle + \langle 1, \bar{\mathcal{L}}_i^r(s,t)\rangle$$
$$\leq \langle 1_{B+\bar{S}_i^r(s,t)}, \bar{\mathcal{Z}}_i^r(s)\rangle + \frac{\varepsilon}{2}.$$

If $I = \varnothing$, then (5.39) implies that $B + \bar{S}_i^r(s,t) \subset B^\varepsilon$. So, (5.45) yields

$$\langle 1_B, \bar{\mathcal{Z}}_i^r(t)\rangle \leq \langle 1_{B^\varepsilon}, \bar{\mathcal{Z}}_i^r(s)\rangle + \frac{\varepsilon}{2}.$$

If $I \neq \varnothing$, then by (5.4), (5.41) and (5.44),

$$\langle 1_B, \bar{\mathcal{Z}}_i^r(t)\rangle \leq \langle 1, \bar{\mathcal{Z}}_i^r(\tau)\rangle + \langle 1, \bar{\mathcal{L}}_i^r(\tau,t)\rangle \leq \langle 1, \bar{\mathcal{Z}}_i^r(\tau)\rangle + \langle 1, \bar{\mathcal{L}}_i^r(s,t)\rangle \leq \varepsilon.$$

In both cases, (5.37) holds. □

5.6. *Tightness.* This section combines the work of Sections 5.1–5.5 to prove the tightness claim of Theorem 4.1.

THEOREM 5.7. *The sequence* $\{(\bar{\mathcal{Z}}^r, \bar{Z}^r, \bar{T}^r, \bar{U}^r, \bar{W}^r)\}$ *is* **C**-*tight.*

PROOF. For each $T > 0$ and $\delta \in [0,T]$, let $\mathbf{w}_T'(\cdot, \delta)$ be the modulus of continuity on $\mathbf{D}([0,\infty), \mathbf{M^I})$ used in Corollary 3.7.4 of [6]. By Definition 5.5,

$$\mathbf{w}_T'(\zeta(\cdot), \delta) \leq \mathbf{w}_{T+\delta}(\zeta(\cdot), 2\delta)$$

for all $\zeta(\cdot) \in \mathbf{D}([0,\infty), \mathbf{M^I})$. So by Lemmas 5.3 and 5.6, the measure valued state descriptors $\{\bar{\mathcal{Z}}^r\}$ satisfy the compact containment and oscillation conditions of Corollary 3.7.4 in [6]. Thus, $\{\bar{\mathcal{Z}}^r\}$ is tight. Moreover, Definition 5.5 and Lemma 5.6 also imply that any weak limit point $\mathcal{Z}$ (obtained as a limit in distribution along a subsequence of $\{\bar{\mathcal{Z}}^r\}$) is continuous almost surely. Since $\bar{Z}^r(\cdot) = \langle 1, \bar{\mathcal{Z}}^r(\cdot)\rangle$ and $\xi \mapsto \langle 1, \xi\rangle$ is continuous on $\mathbf{M}^{\mathbf{I}}$, it follows that $\{(\bar{\mathcal{Z}}^r, \bar{Z}^r)\}$ is **C**-tight.

By (4.3) and (2.1), $\bar{T}^r(\cdot)$ is almost surely Lipschitz continuous with Lipschitz constant $\|C\|$. Since this holds uniformly in $r$, the sequence $\{\bar{T}^r\}$ is tight and any weak limit point $T$ is almost surely Lipschitz continuous with Lipschitz constant $\|C\|$. By (4.4), **C**-tightness of $\{\bar{T}^r\}$ implies **C**-tightness of $\{\bar{U}^r\}$.

As $r \to \infty$, $\bar{W}^r(0) \Rightarrow \langle \chi, \mathcal{Z}^0\rangle$ by (4.14), and $\langle \chi, \bar{\mathcal{L}}^r(\cdot)\rangle \Rightarrow \rho(\cdot)$ by Theorem 5.1. So (4.6) and **C**-tightness of $\{\bar{T}^r\}$ imply **C**-tightness of $\{\bar{W}^r\}$. It follows that $\{(\bar{\mathcal{Z}}^r, \bar{Z}^r, \bar{T}^r, \bar{U}^r, \bar{W}^r)\}$ is **C**-tight. $\square$

5.7. *Weak limits as fluid model solutions.* Let $(\mathcal{Z}, Z, T, U, W)$ be a weak limit of the sequence $\{(\bar{\mathcal{Z}}^r, \bar{Z}^r, \bar{T}^r, \bar{U}^r, \bar{W}^r)\}$, and let $\{q\} \subset \mathcal{R}$ be a subsequence such that

$$(\bar{\mathcal{Z}}^q, \bar{Z}^q, \bar{T}^q, \bar{U}^q, \bar{W}^q) \Rightarrow (\mathcal{Z}, Z, T, U, W) \qquad \text{as } q \to \infty.$$

Note that since $\bar{W}^q(0) = \langle \chi, \bar{\mathcal{Z}}^q(0)\rangle$ for all $q$, assumption (4.14) implies that $(\mathcal{Z}(0), W(0)) \sim (\mathcal{Z}^0, \langle \chi, \mathcal{Z}^0\rangle)$ and so $W(0) \sim \langle \chi, \mathcal{Z}(0)\rangle$. By Theorem 5.1, $(\bar{\mathcal{L}}^q(\cdot), \langle \chi, \bar{\mathcal{L}}^q(\cdot)\rangle \Rightarrow (\nu(\cdot)\vartheta, \rho(\cdot))$ as $q \to \infty$. Using the Skorohod representation theorem, we may assume without loss of generality for the rest of this subsection that $(\mathcal{Z}, Z, T, U, W)$ and $\{(\bar{\mathcal{Z}}^q, \bar{Z}^q, \bar{T}^q, \bar{U}^q, \bar{W}^q, \bar{\mathcal{L}}^q, \langle \chi, \bar{\mathcal{L}}^q\rangle)\}$ are defined on a common probability space $(\Omega, \mathcal{F}, \mathbf{P})$ such that, almost surely, $W(0) = \langle \chi, \mathcal{Z}(0)\rangle$, and as $q \to \infty$,

(5.46) $\qquad (\bar{\mathcal{Z}}^q, \bar{Z}^q, \bar{T}^q, \bar{U}^q, \bar{W}^q, \bar{\mathcal{L}}^q, \langle \chi, \bar{\mathcal{L}}^q(\cdot)\rangle) \to (\mathcal{Z}, Z, T, U, W, \nu(\cdot)\vartheta, \rho(\cdot)),$

uniformly on compact time intervals. Let $\Omega_1$ be the event of probability one on which $W(0) = \langle \chi, \mathcal{Z}(0)\rangle$ and (5.46) holds. For each $q$, let $\Omega_2^q$ be an event of probability one [cf. (5.2) and Lemma 5.2] on which (5.2) holds, and on which for all $i \leq \mathbf{I}$, all $f \in \mathcal{C}_c$, and all finite time intervals $[s, t] \subset [0, \infty)$ satisfying $\inf_{u \in [s,t]} \bar{Z}_i^q(u) > 0$,

$$
(5.47) \quad \begin{aligned}
\langle f, \bar{\mathcal{Z}}_i^q(t)\rangle = \langle f, \bar{\mathcal{Z}}_i^q(s)\rangle &- \int_s^t \langle f', \bar{\mathcal{Z}}_i^q(u)\rangle \frac{\Lambda_i(\bar{Z}^q(u))}{\bar{Z}_i^q(u)} \, du \\
&+ \langle f, \bar{\mathcal{L}}_i^q(t)\rangle - \langle f, \bar{\mathcal{L}}_i^q(s)\rangle.
\end{aligned}
$$

Then $\Omega_2 = \bigcap_q \Omega_2^q$ also has probability one. Define $\Omega_0 = \Omega_1 \cap \Omega_2$.

Lemma 5.8 and Theorem 5.9 below establish that almost surely, $\mathcal{Z}$ is a fluid model solution with auxiliary functions $(Z, T, U, W)$ for the data $(A, C, \Lambda, \nu, \vartheta)$, where $W(t)$ is finite for all $t \geq 0$.

First, recall that for a function $x:[0, \infty) \to \mathbb{R}$, a *regular point* for $x$ is a value of $t \in (0, \infty)$ at which $x$ is differentiable. If $x$ is absolutely continuous, then almost every $t \in (0, \infty)$ is a regular point for $x$, and

$$x(t) = x(0) + \int_0^t \dot{x}(s)\,ds, \qquad t \geq 0,$$

where $\dot{x}$ is equal to the derivative of $x$ whenever $x$ is differentiable, and $\dot{x}$ is equal to zero otherwise. A uniformly Lipschitz continuous function $x : [0, \infty) \to \mathbb{R}$ is absolutely continuous.

LEMMA 5.8. *Almost surely, for all $t \geq 0$, the limit $(\mathcal{Z}, Z, T, U, W)$ satisfies:*

(i) $\|\langle 1_{\{0\}}, \mathcal{Z}(t)\rangle\| = 0$,
(ii) $Z(t) = \langle 1, \mathcal{Z}(t)\rangle$,
(iii) $U(t) = Ct - AT(t)$,
(iv) $W(t) = W(0) + \rho t - T(t)$,
(v) $W(t) = \langle \chi, \mathcal{Z}(t)\rangle$,
(vi) $W$ *is uniformly Lipschitz continuous with Lipschitz constant* $\|\rho\| + \|C\|$,
(vii) *for all $i \leq \mathbf{I}$,*

$$T_i(t) = \int_0^t (\Lambda_i(Z(s))1_{(0,\infty)}(Z_i(s)) + \rho_i 1_{\{0\}}(Z_i(s)))\,ds,$$

(viii) $U_j$ *is nondecreasing for all $j \leq \mathbf{J}$.*

PROOF. Let $T > 0$. It suffices to show (i) for all $t \in [0, T)$. By Lemma 5.4, there exists a sequence $\{a_n : n \in \mathbb{N}\}$ of positive real numbers such that, for each fixed $n$,

$$(5.48) \qquad \liminf_{q \to \infty} \mathbf{P}\left(\sup_{t \in [0,T)} \|\langle 1_{[0,a_n)}, \bar{\mathcal{Z}}^q(t)\rangle\| \leq \frac{1}{n}\right) \geq 1 - \frac{1}{n^2}.$$

For each $n \in \mathbb{N}$, let $\mathbf{A}_n = \{\xi \in \mathbf{M}^{\mathbf{I}} : \|\langle 1_{[0,a_n)}, \xi\rangle\| \leq 1/n\}$, and suppose that $\{\xi^k\} \subset \mathbf{A}_n$ satisfies $\xi^k \xrightarrow{\mathbf{w}} \xi$ as $k \to \infty$. By the Portmanteau theorem,

$$\|\langle 1_{[0,a_n)}, \xi\rangle\| \leq \limsup_{k \to \infty} \|\langle 1_{[0,a_n)}, \xi^k\rangle\| \leq \frac{1}{n}.$$

So $\xi \in \mathbf{A}_n$, which implies that $\mathbf{A}_n \subset \mathbf{M}^{\mathbf{I}}$ is closed for each $n$. By definition of the Skorohod topology, the set

$$\mathbf{B}_n = \{\zeta(\cdot) \in \mathbf{D}([0,\infty), \mathbf{M}^{\mathbf{I}}) : \zeta(t) \in \mathbf{A}_n \text{ for all } t \in [0, T)\}$$

is closed in $\mathbf{D}([0,\infty), \mathbf{M}^{\mathbf{I}})$ for each $n$. Thus, since $\bar{\mathcal{Z}}^q \Rightarrow \mathcal{Z}$, (5.48) and the Portmanteau theorem imply that

$$(5.49) \qquad \mathbf{P}(\mathcal{Z} \in \mathbf{B}_n) \geq \liminf_{q \to \infty} \mathbf{P}(\bar{\mathcal{Z}}^q \in \mathbf{B}_n) \geq 1 - \frac{1}{n^2}.$$

We deduce from (5.49) and the Borel–Cantelli lemma that

$$\mathbf{P}\bigg(\sup_{t\in[0,T)}\|\langle 1_{\{0\}},\mathcal{Z}(t)\rangle\|=0\bigg) \geq \mathbf{P}\bigg(\bigcup_{k=1}^{\infty}\bigcap_{n=k}^{\infty}\{\mathcal{Z}\in\mathbf{B}_n\}\bigg)=1,$$

which proves (i).

Fix an outcome $\omega\in\Omega_0$ and assume for the rest of the proof that all random objects are evaluated at this $\omega$. Property (ii) follows from (4.2) and (5.46). Property (iii) follows from (4.4) and (5.46), and property (iv) follows from (4.6) and (5.46).

To prove (v), fix $t\geq 0$ and $i\leq\mathbf{I}$. Since $\bar{W}_i^q(t)=\langle\chi,\bar{\mathcal{Z}}_i^q(t)\rangle$ for all $q$, and since $\bar{W}_i^q(t)\to W_i(t)$ as $q\to\infty$ by (5.46), it suffices to show that $\langle\chi,\bar{\mathcal{Z}}_i^q(t)\rangle\to\langle\chi,\mathcal{Z}_i(t)\rangle$ as $q\to\infty$. Since $\bar{\mathcal{Z}}_i^q(t)\xrightarrow{\mathbf{w}}\mathcal{Z}_i(t)$ as $q\to\infty$, it suffices to show that the $q$-indexed sequence of measures $\{\bar{\mathcal{Z}}_i^q(t)\}$ is uniformly integrable. To this end, note that if a sequence $\{\xi^q\}\subset\mathbf{M}$ satisfies $\xi^q\xrightarrow{\mathbf{w}}\xi$ and $\langle\chi,\xi^q\rangle\to\langle\chi,\xi\rangle<\infty$, then $\{\xi^q\}$ is uniformly integrable. Thus, $\{\bar{\mathcal{Z}}_i^q(0)\}$ is uniformly integrable by the definition of $\Omega_1$ and (5.46), and $\{\bar{\mathcal{L}}_i^q(t)\}$ is uniformly integrable by (5.46). We conclude from the dynamic equation bound (5.3) that, for each $x>0$,

$$\sup_q\langle\chi 1_{[x,\infty)},\bar{\mathcal{Z}}_i^q(t)\rangle\leq\sup_q(\langle\chi 1_{[x,\infty)},\bar{\mathcal{Z}}_i^q(0)\rangle+\langle\chi 1_{[x,\infty)},\bar{\mathcal{L}}_i^q(t)\rangle).$$

So uniform integrability of $\{\bar{\mathcal{Z}}_i^q(t)\}$ follows from uniform integrability of $\{\bar{\mathcal{Z}}_i^q(0)\}$ and $\{\bar{\mathcal{L}}_i^q(t)\}$.

To prove (vi), recall from the proof of Theorem 5.7 that $T$ is uniformly Lipschitz continuous with Lipschitz constant $\|C\|$. So, (vi) follows from (iv).

To prove (vii), fix $i\leq\mathbf{I}$. Since $W_i$ and $T_i$ are uniformly Lipschitz continuous, they are both absolutely continuous. Let $t>0$ be a regular point for both $W_i$ and $T_i$. Then $\dot{W}_i(t)=\rho_i-\dot{T}_i(t)$ by (iv). If $Z_i(t)=0$, then $W_i(t)=0$ by (v). Since $W_i$ is a nonnegative function, this implies that $\dot{W}_i(t)=0$ and so $\dot{T}_i(t)=\rho_i$. Alternatively, suppose that $Z_i(t)>0$. By (ii), continuity of $\mathcal{Z}_i$ implies continuity of $Z_i$. So $Z_i(s)>0$ for all $s\in[t,t+h]$ and all sufficiently small $h>0$. In this case, (5.46), (4.3), continuity of $\Lambda_i$ on $\{z\in\mathbb{R}_+^{\mathbf{I}}:z_i>0\}$, (2.1), and the bounded convergence theorem, imply that

$$T_i(t+h)-T_i(t)=\lim_{q\to\infty}(\bar{T}_i^q(t+h)-\bar{T}_i^q(t))$$

(5.50)
$$=\lim_{q\to\infty}\int_t^{t+h}\Lambda_i(\bar{Z}^q(s))\,ds$$

$$=\int_t^{t+h}\Lambda_i(Z(s))\,ds.$$

Since $\Lambda_i(Z(\cdot))$ is continuous at $t$ for $Z_i(t)>0$, it follows that

(5.51) $$\dot{T}_i(t)=\begin{cases}\Lambda_i(Z(t)),&\text{if }Z_i(t)>0,\\\rho_i,&\text{if }Z_i(t)=0.\end{cases}$$

Since almost every $t > 0$ is a regular point for $W_i$ and $T_i$, (5.51) implies (vii).

Property (viii) follows because $\bar{U}_j^q$ is nondecreasing for each $q$ and $j \leq \mathbf{J}$, and because $\bar{U}^q \to U$ uniformly on compact time intervals by (5.46). □

The next result establishes the family of dynamic equations satisfied by the limit in (5.46).

THEOREM 5.9. *Almost surely, for all $i \leq \mathbf{I}$, $f \in \mathcal{C}$, and $t \geq 0$,*

$$\langle f, \mathcal{Z}_i(t)\rangle = \langle f, \mathcal{Z}_i(0)\rangle - \int_0^t \langle f', \mathcal{Z}_i(s)\rangle \frac{\Lambda_i(Z(s))}{Z_i(s)} \, ds$$
(5.52)
$$+ \nu_i \langle f, \vartheta_i\rangle \int_0^t 1_{(0,\infty)}(Z_i(s)) \, ds.$$

Recall that the first integrand above is defined to be zero when $Z_i(s) = 0$.

PROOF OF THEOREM 5.9. All random objects in this proof are evaluated at a fixed outcome $\omega \in \Omega_0$ such that (5.2), (5.46), (5.47) and the properties listed in Lemma 5.8 hold. The theorem will be proved first for $f \in \mathcal{C}_c$, and an extension to $\mathcal{C}$ is made at the end. Recall that for $f \in \mathcal{C}_c$, the derivative $f'$ has compact support and $\|f'\|_\infty < \infty$. Note also that since $f(0) = f'(0) = 0$, there exists a constant $C_f < \infty$ such that $|f(x)| \leq C_f x$ for all $x \in \mathbb{R}_+$. Therefore,

(5.53) $\quad |\langle f, \mathcal{Z}_i(t)\rangle| \leq C_f \langle \chi, \mathcal{Z}_i(t)\rangle = C_f W_i(t) \quad$ for all $t \geq 0$.

The following preliminary result is used several times in this proof.

For each fixed $f \in \mathcal{C}_c$, each $i \leq \mathbf{I}$ and each interval $[s,t] \subset \mathbb{R}_+$ satisfying $\inf_{u \in [s,t]} Z_i(u) > 0$, we have

(5.54) $\quad \langle f, \mathcal{Z}_i(t)\rangle = \langle f, \mathcal{Z}_i(s)\rangle - \int_s^t \langle f', \mathcal{Z}_i(u)\rangle \frac{\Lambda_i(Z(u))}{Z_i(u)} \, du + \nu_i(t-s)\langle f, \vartheta_i\rangle.$

To see this, fix an $f$, $i$ and interval $[s,t]$ satisfying the assumption. By (5.46) and Lemma 5.8 (ii), $\bar{Z}_i^q(\cdot) \to Z_i(\cdot)$ as $q \to \infty$, uniformly on compact time intervals. Thus,

(5.55) $\qquad\qquad\qquad \liminf_{q \to \infty} \inf_{u \in [s,t]} \bar{Z}_i^q(u) > 0.$

By definition of $\Omega_2$, this implies that (5.47) holds for all sufficiently large $q$. Let $q \to \infty$ in (5.47). Note that $f$ and $f'$ are elements of $\mathbf{C}_b(\mathbb{R}_+)$, and that $\Lambda_i$ is continuous on $\{z \in \mathbb{R}_+^\mathbf{I} : z_i > 0\}$ (see Definition 2.1). Thus, (5.46) and (5.55) imply that for all $q$ sufficiently large, the integrands in the second right-hand term of (5.47) are uniformly bounded on $[s,t]$ and converge pointwise on $[s,t]$ to $\langle f', \mathcal{Z}_i(\cdot)\rangle \frac{\Lambda_i(Z(\cdot))}{Z_i(\cdot)}$. Apply bounded convergence to this term and apply (5.46) to the remaining terms in (5.47) to obtain (5.54).

We now proceed with the proof of the theorem. Fix $f \in \mathcal{C}_c$, $i \leq \mathbf{I}$, and an interval $[s,t] \subset \mathbb{R}_+$ with $t > s$. Define

$$\tau_0 = \inf\{u \in [s,t] : Z_i(u) = 0\}, \tag{5.56}$$

where the infimum of the empty set is defined to be $t$. If $\tau_0 > s$, then $\inf_{u \in [s,\tau]} Z_i(u) > 0$ for all intervals $[s,\tau] \subset [s,\tau_0)$. So, for each such $\tau$, (5.54) with $\tau$ in place of $t$ implies that

$$|\langle f, \mathcal{Z}_i(\tau)\rangle - \langle f, \mathcal{Z}_i(s)\rangle| \leq (\tau - s)\|f'\|_\infty \|C\| + (\tau - s)\nu_i \|f\|_\infty. \tag{5.57}$$

Since both sides of (5.57) are continuous in $\tau$, letting $\tau \uparrow \tau_0$ yields

$$|\langle f, \mathcal{Z}_i(\tau_0)\rangle - \langle f, \mathcal{Z}_i(s)\rangle| \leq (\tau_0 - s)\|f'\|_\infty \|C\| + (\tau_0 - s)\nu_i \|f\|_\infty. \tag{5.58}$$

Note that if $\tau_0 = s$, then (5.58) holds trivially. If $\tau_0 < t$, then $Z_i(\tau_0) = 0$ by continuity of $Z_i$, and so $\langle f, \mathcal{Z}_i(\tau_0)\rangle = W_i(\tau_0) = 0$. Then by (5.53) and property (vi) of Lemma 5.8,

$$\begin{aligned}
|\langle f, \mathcal{Z}_i(t)\rangle - \langle f, \mathcal{Z}_i(\tau_0)\rangle| = |\langle f, \mathcal{Z}_i(t)\rangle| &\leq C_f W_i(t) \\
&= C_f |W_i(t) - W_i(\tau_0)| \\
&\leq C_f(\|\rho\| + \|C\|)(t - \tau_0).
\end{aligned} \tag{5.59}$$

If $\tau_0 = t$, then (5.59) holds trivially. Combining (5.58) and (5.59) yields

$$\begin{aligned}
|\langle f, \mathcal{Z}_i(t)\rangle - \langle f, \mathcal{Z}_i(s)\rangle| &\leq |\langle f, \mathcal{Z}_i(t)\rangle - \langle f, \mathcal{Z}_i(\tau_0)\rangle| \\
&\quad + |\langle f, \mathcal{Z}_i(\tau_0)\rangle - \langle f, \mathcal{Z}_i(s)\rangle| \\
&\leq C_f(\|\rho\| + \|C\|)(t - \tau_0) \\
&\quad + (\|f'\|_\infty \|C\| + \nu_i \|f\|_\infty)(\tau_0 - s) \\
&\leq K_f(t - s),
\end{aligned}$$

where

$$K_f = C_f(\|\rho\| + \|C\|) + \|f'\|_\infty \|C\| + \nu_i \|f\|_\infty.$$

Since $s \leq t$ were arbitrary, it follows that $\langle f, \mathcal{Z}_i(\cdot)\rangle$ is uniformly Lipschitz continuous and is therefore absolutely continuous on $\mathbb{R}_+$. Suppose $t > 0$ is a regular point for both $\langle f, \mathcal{Z}_i(\cdot)\rangle$ and $W_i(\cdot)$. If $Z_i(t) > 0$, then by (5.54) (with $[t, t+h]$ in place of $[s,t]$), and continuity of $\Lambda_i$ on $\{z \in \mathbb{R}_+^\mathbf{I} : z_i > 0\}$ (see Definition 2.1),

$$\lim_{h \to 0} \frac{\langle f, \mathcal{Z}_i(t+h)\rangle - \langle f, \mathcal{Z}_i(t)\rangle}{h} = -\langle f', \mathcal{Z}_i(t)\rangle \frac{\Lambda_i(Z(t))}{Z_i(t)} + \nu_i \langle f, \vartheta_i\rangle. \tag{5.60}$$

If $Z_i(t) = 0$, then $\langle f, \mathcal{Z}_i(t) \rangle = W_i(t) = 0$, and so $\dot{W}_i(0) = 0$ because $W_i$ is a nonnegative function. Combining this with (5.53), we obtain

$$
\begin{aligned}
\limsup_{h \to 0} \left| \frac{\langle f, \mathcal{Z}_i(t+h) \rangle - \langle f, \mathcal{Z}_i(t) \rangle}{h} \right| &= \limsup_{h \to 0} \left| \frac{\langle f, \mathcal{Z}_i(t+h) \rangle}{h} \right| \\
&\leq \limsup_{h \to 0} C_f \left| \frac{W_i(t+h)}{h} \right| \\
&= \limsup_{h \to 0} C_f \left| \frac{W_i(t+h) - W_i(t)}{h} \right| \\
&= C_f |\dot{W}_i(t)| = 0.
\end{aligned}
\tag{5.61}
$$

Combining (5.60) and (5.61), we obtain

$$
\frac{d}{dt}\langle f, \mathcal{Z}_i(t) \rangle = \begin{cases} -\langle f', \mathcal{Z}_i(t) \rangle \dfrac{\Lambda_i(Z(t))}{Z_i(t)} + \nu_i \langle f, \vartheta_i \rangle, & Z_i(t) > 0, \\ 0, & Z_i(t) = 0. \end{cases}
\tag{5.62}
$$

The set of all $t \in (0, \infty)$ that are regular points for both $\langle f, \mathcal{Z}_i(\cdot) \rangle$ and $W_i(\cdot)$ has full Lebesgue measure, so (5.52) follows from (5.62).

This proves the theorem for $f \in \mathcal{C}_c$. To extend to $\mathcal{C}$, choose functions $\{g_n : n \in \mathbb{N}\} \subset \mathbf{C}_b^1(\mathbb{R}_+)$ such that $1_{[0,n]} \leq g_n \leq 1_{[0,n+1]}$ and $\|g'_n\|_\infty \leq 2$ for all $n$. For $f \in \mathcal{C}$, define $f_n = f g_n$ so that $f_n \in \mathcal{C}_c$ for all $n$. Then for all $n \in \mathbb{N}$ and $t \geq 0$,

$$
\begin{aligned}
\langle f_n, \mathcal{Z}_i(t) \rangle &= \langle f_n, \mathcal{Z}_i(0) \rangle - \int_0^t \langle f'_n, \mathcal{Z}_i(s) \rangle \frac{\Lambda_i(Z(s))}{Z_i(s)} \, ds \\
&\quad + \nu_i \langle f_n, \vartheta_i \rangle \int_0^t 1_{(0,\infty)}(Z_i(s)) \, ds.
\end{aligned}
\tag{5.63}
$$

Since $f_n \to f$ pointwise and boundedly as $n \to \infty$, the bounded convergence theorem implies that the left-hand side, as well as the first and third terms on the right-hand side of (5.63), converge to the corresponding terms of (5.52) as $n \to \infty$. Similarly, $f'_n \to f'$ pointwise and boundedly as $n \to \infty$. So, the integrand in the second right-hand term of (5.63) converges pointwise on $[0,t]$ to $\langle f', \mathcal{Z}_i(\cdot) \rangle \frac{\Lambda_i(Z(\cdot))}{Z_i(\cdot)}$ as $n \to \infty$. For each $n \in \mathbb{N}$,

$$
\sup_{s \in [0,t]} \left| \langle f'_n, \mathcal{Z}_i(s) \rangle \frac{\Lambda_i(Z(s))}{Z_i(s)} \right| \leq \|f'_n\|_\infty \|C\| \leq (\|f'\|_\infty + 2\|f\|_\infty) \|C\|.
$$

So, the bounded convergence theorem implies that the second right-hand term in (5.63) converges to the corresponding term in (5.52) as $n \to \infty$. □

PROOF OF THEOREM 4.1. The sequence $\{(\bar{\mathcal{Z}}^r, \bar{Z}^r, \bar{T}^r, \bar{U}^r, \bar{W}^r)\}$ is **C**-tight by Theorem 5.7. Let $(\mathcal{Z}, Z, T, U, W)$ be a weak limit point of this sequence. Then by Theorem 5.9 and properties (i)–(v), (vii) and (viii) of Lemma 5.8, $\mathcal{Z}$ is almost surely a fluid model solution with auxiliary functions $(Z, T, U, W)$ for the data $(A, C, \Lambda, \nu, \vartheta)$ (see Definitions 3.1 and 3.2), and $W(t)$ is finite for all $t \geq 0$. $\square$

**6. Invariant states for fluid model under weighted $\alpha$-fair policies.** In this section, we consider the special case of weighted $\alpha$-fair policies. Fix fluid model data $(A, C, \Lambda, \nu, \vartheta)$, where $\Lambda$ is a weighted $\alpha$-fair bandwidth sharing policy with parameters $(\alpha, \kappa)$ as described in the example of Section 2.2. Under a natural condition on the network parameters $A$, $C$, $\nu$ and $\vartheta$, there exist fluid model solutions that are time invariant. This section identifies the condition and characterizes the set of these *invariant states*.

The following representation of the weighted $\alpha$-fair policy $\Lambda$ follows from Lemma A.4 of [12]. (Although it is assumed at the beginning of [12] that $A$ has full row rank, this property is not used in the proof of Lemma A.4 in [12], and hence the result holds without restriction on the row rank of $A$.)

PROPOSITION 6.1. *For each $z \in \mathbb{R}_+^{\mathbf{I}}$, there exists at least one $p \in \mathbb{R}_+^{\mathbf{J}}$, depending on $z$, such that*

$$(6.1) \qquad \Lambda_i(z) = z_i \left( \frac{\kappa_i}{\sum_{j \leq \mathbf{J}} p_j A_{ji}} \right)^{1/\alpha} \qquad \text{for all } i \in \mathcal{I}_+(z),$$

*where*

$$(6.2) \qquad p_j \left( C_j - \sum_{i \leq \mathbf{I}} A_{ji} \Lambda_i(z) \right) = 0 \qquad \text{for all } j \leq \mathbf{J}.$$

The $(p_j : j \leq \mathbf{J})$ are Lagrange multipliers for the optimization problem (2.3)–(2.5), one for each of the capacity constraints in (2.4). Note that for each $i \in \mathcal{I}_+(z)$, the bandwidth $\Lambda_i(z) > 0$ by Definition 2.1(i), and $z_i > 0$ by definition. Thus, (6.1) implies that the denominator on the right-hand side of (6.1) does not vanish.

DEFINITION 6.2. A vector of measures $\xi \in \mathbf{M}^{\mathbf{I}}$ is an invariant state for the fluid model if there is a fluid model solution $\zeta(\cdot)$ satisfying $\zeta(t) = \xi$ for all $t \geq 0$.

The following notation helps describe invariant states. Recall that $\mu_i = \langle \chi, \vartheta_i \rangle^{-1}$ and $\rho_i = \nu_i / \mu_i$ for $i \leq \mathbf{I}$. Also recall that $\mathcal{I}_+(z) = \{i \leq \mathbf{I} : z_i > 0\}$ and $\mathcal{I}_0(z) = \{i \leq \mathbf{I} : z_i = 0\}$ for $z \in \mathbb{R}_+^{\mathbf{I}}$. Let

$$\mathcal{P} = \{z \in \mathbb{R}_+^{\mathbf{I}} : \Lambda_i(z) = \rho_i \text{ for all } i \in \mathcal{I}_+(z)\}.$$

For each $i \leq \mathbf{I}$, let $\vartheta_i^e$ denote the excess lifetime distribution associated with $\vartheta_i$. The probability measure $\vartheta_i^e$ is absolutely continuous with density

$$p_i^e(x) = \mu_i \langle 1_{(x,\infty)}, \vartheta_i \rangle, \qquad x \in \mathbb{R}_+.$$

THEOREM 6.3. *There exist invariant states for the fluid model if and only if*

(6.3) $$A\rho \leq C.$$

*When* (6.3) *holds, the set of invariant states is given by*

(6.4) $$\mathcal{M} = \{\xi \in \mathbf{M}^{\mathbf{I}} : \xi_i = z_i \vartheta_i^e \text{ for all } i \leq \mathbf{I} \text{ and some } z \in \mathcal{P}\}.$$

PROOF. Suppose that $\xi$ is an invariant state and let $\zeta(\cdot) \equiv \xi$ be the corresponding fluid model solution with auxiliary functions $(z, \tau, u)$ given by Definition 3.1 (we omit $w$ here). Then $z$ is a constant vector, denoted $z = \langle 1, \xi \rangle$. For each $f \in \mathcal{C}$, $i \in \mathcal{I}_+(z)$ and $t \geq 0$, property (iii) of Definition 3.2 yields

(6.5) $$\langle f, \xi_i \rangle = \langle f, \xi_i \rangle - t \langle f', \xi_i \rangle \frac{\Lambda_i(z)}{z_i} + t\nu_i \langle f, \vartheta_i \rangle.$$

Since (cf. Proposition 3.1 in [20]),

(6.6) $$\langle f, \vartheta_i \rangle = \frac{1}{\mu_i} \langle f', \vartheta_i^e \rangle \qquad \text{for all } f \in \mathcal{C},$$

canceling like terms in (6.5) yields

(6.7) $$\langle f', \xi_i \rangle \frac{\Lambda_i(z)}{z_i} = \rho_i \langle f', \vartheta_i^e \rangle.$$

Replacing $f$ by a suitable sequence of functions $\{f_n\} \subset \mathcal{C}$ satisfying $f_n' \geq 0$ and $f_n' \uparrow 1_{(0,\infty)}$, the monotone convergence theorem implies that (6.7) holds with $f' \equiv 1_{(0,\infty)}$. So by property (i) of Definition 3.2, and since $\vartheta_i^e$ does not charge $\{0\}$,

(6.8) $$\Lambda_i(z) = \rho_i \qquad \text{for all } i \in \mathcal{I}_+(z).$$

Since $z_i = 0$ for $i \in \mathcal{I}_0(z)$, this implies by Definition 3.1 that for all $t \geq 0$,

(6.9) $$\tau_i(t) = \int_0^t (\Lambda_i(z) 1_{(0,\infty)}(z_i) + \rho_i 1_{\{0\}}(z_i)) \, ds = \rho_i t \qquad \text{for all } i \leq \mathbf{I}.$$

Thus, $u(t) = Ct - A\rho t = (C - A\rho)t$ for all $t \geq 0$. Since $u$ is nondecreasing by property (ii) of Definition 3.2, (6.3) holds. Moreover, substituting (6.8) into (6.7) and canceling $\rho_i$ yields

$$\frac{\langle f', \xi_i \rangle}{z_i} = \langle f', \vartheta_i^e \rangle \qquad \text{for all } f \in \mathcal{C} \text{ and } i \in \mathcal{I}_+(z).$$

This implies (in the same manner as in the proof of Theorem 1.1 in [20]) that

(6.10) $$\xi_i = z_i \vartheta_i^e \qquad \text{for all } i \in \mathcal{I}_+(z).$$

Since $\xi_i = \mathbf{0}$ for all $i \in \mathcal{I}_0(z)$, combining (6.10) with (6.8) implies that $\xi \in \mathcal{M}$.

To prove the converse, suppose that (6.3) holds and let $\xi \in \mathcal{M}$. Define $\zeta(t) = \xi$ for all $t \geq 0$, and let $(z, \tau, u)$ be the auxiliary functions of $\zeta$ given by Definition 3.1. Since $\vartheta_i^e$ does not charge $\{0\} \subset \mathbb{R}_+$ for each $i \leq \mathbf{I}$, $\zeta(\cdot)$ satisfies property (i) of Definition 3.2. Note that $z = \langle 1, \xi \rangle$ is a constant vector, and that $z \in \mathcal{P}$ because $\xi \in \mathcal{M}$. Thus, $\tau$ satisfies (6.9) and so $u(t) = Ct - A\rho t = (C - A\rho)t$ for all $t \geq 0$. By (6.3), $u_j$ is nondecreasing for all $j \leq \mathbf{J}$, and so property (ii) of Definition 3.2 holds. Let $f \in \mathcal{C}$ and $i \in \mathcal{I}_+(z)$. Since $\xi \in \mathcal{M}$ and $z \in \mathcal{P}$, (6.6) implies that for all $t \geq 0$,

$$t\langle f', \xi_i \rangle \frac{\Lambda_i(z)}{z_i} = t\rho_i \langle f', \vartheta_i^e \rangle = t\nu_i \langle f, \vartheta_i \rangle.$$

Thus, (6.5) holds and so $\zeta$ satisfies property (iii) of Definition 3.2 for $i \in \mathcal{I}_+(z)$. This property holds for $i \notin \mathcal{I}_+(z)$ since then all terms are zero. Thus, $\zeta$ is a fluid model solution. Note that $\mathcal{M}$ is nonempty because $\xi = \mathbf{0}$ is in $\mathcal{M}$. □

Under condition (6.3), the set $\mathcal{P}$ can be characterized using results in Kelly and Williams [12]. (Although it was assumed in [12] that $A$ has full row rank, as we explain below, the results that we cite below from [12] hold without that additional restriction.) Let

$$\mathcal{J}_* = \left\{ j \leq \mathbf{J} : \sum_{i \leq \mathbf{I}} A_{ji} \rho_i = C_j \right\}$$

and let $\mathbf{J}_* = |\mathcal{J}_*|$. For $z \in \mathbb{R}_+^\mathbf{I}$, define

$$F(z) = \frac{1}{\alpha + 1} \sum_{i \leq \mathbf{I}} \nu_i \kappa_i \mu_i^{\alpha - 1} \left( \frac{z_i}{\nu_i} \right)^{\alpha + 1}.$$

When $\mathcal{J}_* \neq \varnothing$, define $\Delta : \mathbb{R}_+^{\mathbf{J}_*} \to \mathbb{R}_+^\mathbf{I}$ by

$$\Delta(w) = \arg\min \left\{ F(z) : z \in \mathbb{R}_+^\mathbf{I} \text{ and } \sum_{i \leq \mathbf{I}} A_{ji} \frac{z_i}{\mu_i} \geq w_j \text{ for all } j \in \mathcal{J}_* \right\}.$$

For each $j \in \mathcal{J}_*$, $A_{ji} > 0$ for some $i \leq \mathbf{I}$. It follows that the feasible set for the above optimization problem is nonempty. Then since $F$ is nonnegative, continuous, strictly convex and satisfies $F(z) \to \infty$ as $\|z\| \to \infty$, $\Delta(w)$ is well defined as the unique minimum of the optimization problem.

LEMMA 6.4. *Assume that* (6.3) *holds. If* $\mathcal{J}_* = \varnothing$, *then* $\mathcal{P} = \{0\}$ *and the only invariant state is* $\xi = \mathbf{0}$. *If* $\mathcal{J}_* \neq \varnothing$, *then the following three conditions are equivalent:*

(i) $z \in \mathcal{P}$,
(ii) *for some* $q \in \mathbb{R}_+^{\mathbf{J}_*}$, $z_i = \rho_i(\frac{1}{\kappa_i}\sum_{j \in \mathcal{J}_*} q_j A_{ji})^{1/\alpha}$ *for all* $i \leq \mathbf{I}$,
(iii) $z = \Delta(w(z))$, *where* $w_j(z) = \sum_{i \leq \mathbf{I}} A_{ji}\frac{z_i}{\mu_i}$ *for all* $j \in \mathcal{J}_*$.

PROOF. If $\mathcal{J}_* = \varnothing$, then $A\rho < C$. In this case, there is no $z \neq 0$ such that $\Lambda_i(z) = \rho_i$ for all $i \in \mathcal{I}_+(z)$. This follows because for $z \neq 0$, the optimal solution $\Lambda(z)$ of the concave optimization problem (2.3)–(2.5) must have one of the constraints binding, that is, $(A\Lambda(z))_j = C_j$ for some $j \leq \mathbf{J}$. (For this, we use the fact that there is at least one route $i$, which necessarily uses at least one resource $j$, and hence there is at least one $i$ and $j$ such that $A_{ji} > 0$.) It follows that $\mathcal{P} = \{0\}$ when $\mathcal{J}_* = \varnothing$.

If $\mathcal{J}_* \neq \varnothing$, then the three equivalent characterizations of $\mathcal{P}$ follow from Lemma 5.1 and Theorems 5.1 and 5.3 in [12]. (At the beginning of [12], it is assumed that $A$ has full row rank. However, scrutiny of the proofs of Lemma 5.1 and Theorems 5.1 and 5.3 of that paper reveals that the above equivalence still holds without this additional assumption. Indeed, one only needs the weaker property that for each $j \in \mathcal{J}_*$, the $j$th row of $A$ has at least one nonzero entry. This property holds here by the definition of $\mathcal{J}_*$.) □

## APPENDIX

PROOF OF THEOREM 5.1. For each $r \in \mathcal{R}$, define a simplified fluid scaled load process

$$\bar{\mathcal{V}}_i^r(t) = \frac{1}{r}\sum_{k=1}^{\lfloor rt \rfloor} \delta_{v_{ik}^r}, \qquad t \geq 0, i \leq \mathbf{I} \tag{A.1}$$

and let $\vartheta(t) = t\vartheta$ for all $t \geq 0$. We first show that

$$(\bar{\mathcal{V}}^r(\cdot), \langle \chi, \bar{\mathcal{V}}^r(\cdot) \rangle) \Rightarrow (\vartheta(\cdot), \langle \chi, \vartheta(\cdot) \rangle) \qquad \text{as } r \to \infty. \tag{A.2}$$

Since for each $r \in \mathcal{R}$, $t \geq 0$, and $i \leq \mathbf{I}$, $\langle \chi, \vartheta_i \rangle = 1/\mu_i$ and

$$\langle \chi, \bar{\mathcal{V}}_i^r(t) \rangle = \frac{\lfloor rt \rfloor}{r}\frac{1}{\lfloor rt \rfloor}\sum_{k=1}^{\lfloor rt \rfloor}\left(v_{ik}^r - \frac{1}{\mu_i^r}\right) + \frac{\lfloor rt \rfloor}{r\mu_i^r},$$

the second component of (A.2) follows from assumptions (4.12), (4.13) and a functional weak law of large numbers. Note that for each $K > 0$, the set $\{\xi \in \mathbf{M}^{\mathbf{I}}: \|\langle 1 \vee \chi, \xi \rangle\| \leq K\}$ is relatively compact in $\mathbf{M}^{\mathbf{I}}$ (see [11], Theorem

15.7.5). For each $T > 0$, $\sup_{t \in [0,T]} \|\langle 1 \vee \chi, \bar{\mathcal{V}}^r(t) \rangle\| \leq T \vee \|\langle \chi, \bar{\mathcal{V}}^r(T) \rangle\|$, and so the second component of (A.2) implies the compact containment condition

(A.3) $$\lim_{K \to \infty} \liminf_{r \to \infty} \mathbf{P}^r \left( \sup_{t \in [0,T]} \|\langle 1 \vee \chi, \bar{\mathcal{V}}^r(t) \rangle\| \leq K \right) = 1.$$

Moreover, for all $r \in \mathcal{R}$, $i \leq \mathbf{I}$, $t \geq s \geq 0$, and all nonempty closed $B \subset \mathbb{R}_+$, (A.1) implies the two inequalities

$$\langle 1_B, \bar{\mathcal{V}}_i^r(s) \rangle \leq \langle 1_B, \bar{\mathcal{V}}_i^r(t) \rangle \leq \langle 1_{B^{t-s}}, \bar{\mathcal{V}}_i^r(t) \rangle + t - s,$$

$$\langle 1_B, \bar{\mathcal{V}}_i^r(t) \rangle = \langle 1_B, \bar{\mathcal{V}}_i^r(s) \rangle + \frac{1}{r} \sum_{k=\lfloor rs \rfloor + 1}^{\lfloor rt \rfloor} 1_B(v_{ik}^r) \leq \langle 1_{B^{t-s}}, \bar{\mathcal{V}}_i^r(s) \rangle + t - s.$$

So by (1.1) and (1.2),
$$\mathbf{d_I}[\bar{\mathcal{V}}^r(s), \bar{\mathcal{V}}^r(t)] \leq t - s \qquad \text{for all } r \in \mathcal{R} \text{ and } t \geq s \geq 0.$$

On combining this with (A.3), we see that $\{\bar{\mathcal{V}}^r(\cdot)\}$ is $\mathbf{C}$-tight. Let $\{\bar{\mathcal{V}}^q(\cdot)\} \subset \{\bar{\mathcal{V}}^r(\cdot)\}$ be a weakly convergent subsequence with almost surely continuous limit $\mathcal{V}(\cdot)$. Then by the continuous mapping theorem, for all $f \in \mathbf{C}_b(\mathbb{R}_+)$,

(A.4) $$\langle f, \bar{\mathcal{V}}^q(\cdot) \rangle \Rightarrow \langle f, \mathcal{V}(\cdot) \rangle \qquad \text{as } q \to \infty.$$

On the other hand, for each $q$, $f \in \mathbf{C}_b(\mathbb{R}_+)$, $t \geq 0$, and $i \leq \mathbf{I}$,

$$\langle f, \bar{\mathcal{V}}_i^q(t) \rangle = \frac{\lfloor qt \rfloor}{q} \frac{1}{\lfloor qt \rfloor} \sum_{k=1}^{\lfloor qt \rfloor} (f(v_{ik}^q) - \langle f, \vartheta_i^q \rangle) + \frac{\lfloor qt \rfloor}{q} \langle f, \vartheta_i^q \rangle.$$

Assumptions (4.11)–(4.13), and a functional weak law of large numbers imply that $\langle f, \bar{\mathcal{V}}^q(\cdot) \rangle \Rightarrow \langle f, \vartheta(\cdot) \rangle$ as $q \to \infty$ for each $f \in \mathbf{C}_b(\mathbb{R}_+)$. Combining with (A.4), we see that $\mathcal{V}(\cdot) \equiv \vartheta(\cdot)$ almost surely, and so $\bar{\mathcal{V}}^r(\cdot) \Rightarrow \vartheta(\cdot)$ as $r \to \infty$. Since the limits are deterministic, the convergence in (A.2) is indeed joint.

By assumption (4.10), $\bar{E}^r(\cdot) \Rightarrow \nu(\cdot)$ as $r \to \infty$. Since $\bar{\mathcal{L}}^r(\cdot) = \bar{\mathcal{V}}^r(\bar{E}^r(\cdot))$ and $\langle \chi, \vartheta(\nu(\cdot)) \rangle = \rho(\cdot)$, (A.2) and the random time change theorem imply (5.1). □

PROOF OF LEMMA 5.2. Let $\Omega_1^r$ be an event of probability one on which (5.2) holds and fix $\omega \in \Omega_1^r$. For the rest of the proof, all random objects are evaluated at this particular $\omega$. Fix $i \leq \mathbf{I}$, $f \in \mathcal{C}_c$, and let $[s, t]$ be an interval satisfying $\inf_{u \in [s,t]} \bar{Z}_i^r(u) > 0$. It suffices to show (5.7). Since $\bar{Z}_i^r(\cdot)$ is right continuous with finite left limits, there exist $\varepsilon, M \in (0, \infty)$ such that

(A.5) $$\varepsilon \leq \inf_{u \in [s,t]} \bar{Z}_i^r(u) \leq \sup_{u \in [s,t]} \bar{Z}_i^r(u) \leq M.$$

Let $l = t - s$ and, for $n, j \in \mathbb{N}$, let $t_j = s + jl/n$ and $t^j = t_{j+1}$. For each $n$,

$$\langle f, \bar{\mathcal{Z}}_i^r(t) \rangle - \langle f, \bar{\mathcal{Z}}_i^r(s) \rangle = \sum_{j=0}^{n-1} (\langle f, \bar{\mathcal{Z}}_i^r(t^j) \rangle - \langle f, \bar{\mathcal{Z}}_i^r(t_j) \rangle).$$

Add and subtract a term in each summand to get

$$\langle f, \bar{\mathcal{Z}}_i^r(t)\rangle - \langle f, \bar{\mathcal{Z}}_i^r(s)\rangle = \sum_{j=0}^{n-1}(\langle f, \bar{\mathcal{Z}}_i^r(t^j)\rangle - \langle f(\cdot - \bar{S}_i^r(t_j, t^j)), \bar{\mathcal{Z}}_i^r(t_j)\rangle)$$
$$+ \sum_{j=0}^{n-1}(\langle f(\cdot - \bar{S}_i^r(t_j, t^j)), \bar{\mathcal{Z}}_i^r(t_j)\rangle - \langle f, \bar{\mathcal{Z}}_i^r(t_j)\rangle).$$

Use the dynamic equation (5.2) in the first term and rewrite the second term on the right to obtain

$$\langle f, \bar{\mathcal{Z}}_i^r(t)\rangle - \langle f, \bar{\mathcal{Z}}_i^r(s)\rangle = \sum_{j=0}^{n-1} \frac{1}{r} \sum_{k=r\bar{E}_i^r(t_j)+1}^{r\bar{E}_i^r(t^j)} f(v_{ik}^r - \bar{S}_i^r(U_{ik}^r r^{-1}, t^j))$$

(A.6)

$$+ \sum_{j=0}^{n-1} \langle f(\cdot - \bar{S}_i^r(t_j, t^j)) - f(\cdot), \bar{\mathcal{Z}}_i^r(t_j)\rangle.$$

Denote the first and second right-hand terms in (A.6) by $a_n^r$ and $b_n^r$, respectively, and consider first $a_n^r$. Since $f \in \mathcal{C}_c$, a first-order Taylor expansion of each summand yields

(A.7) $$f(v_{ik}^r - \bar{S}_i^r(U_{ik}^r r^{-1}, t^j)) = f(v_{ik}^r) + f'(w_j^k)h_j^k,$$

where for each $j$ and $k$, $h_j^k = -\bar{S}_i^r(U_{ik}^r r^{-1}, t^j)$ and $w_j^k \in \mathbb{R}$ is in the interval $[v_{ik}^r + h_j^k, v_{ik}^r]$. Since $U_{ik}^r r^{-1} \in (t_j, t^j]$ for each pair $j, k$ in (A.7), (2.1) and (A.5) imply that

(A.8) $$\max_{j,k} |h_j^k| \leq \max_j \int_{t_j}^{t^j} \frac{\Lambda_i(\bar{Z}^r(u))}{\bar{Z}_i^r(u)} du \leq \frac{l\|C\|}{n\varepsilon}.$$

Using (A.7) and (A.8), deduce that for each $n$,

$$|a_n^r - (\langle f, \bar{\mathcal{L}}_i^r(t)\rangle - \langle f, \bar{\mathcal{L}}_i^r(s)\rangle)| = \left| a_n^r - \frac{1}{r} \sum_{k=r\bar{E}_i^r(s)+1}^{r\bar{E}_i^r(t)} f(v_{ik}^r) \right|$$

$$= \left| \sum_{j=0}^{n-1} \frac{1}{r} \sum_{k=r\bar{E}_i^r(t_j)+1}^{r\bar{E}_i^r(t^j)} f'(w_j^k) h_j^k \right|$$

$$\leq (\bar{E}_i^r(t) - \bar{E}_i^r(s))\|f'\|_\infty \frac{l\|C\|}{n\varepsilon}.$$

So as $n \to \infty$,

(A.9) $$a_n^r \to (\langle f, \bar{\mathcal{L}}_i^r(t)\rangle - \langle f, \bar{\mathcal{L}}_i^r(s)\rangle).$$

Next, consider $b_n^r$. Another first-order Taylor expansion for each $x \in \mathbb{R}_+$ and $j \in \{0, \ldots, n-1\}$ yields

$$f(x - \bar{S}_i^r(t_j, t^j)) - f(x) = f'(w_j^x)h_j, \tag{A.10}$$

where $h_j = -\bar{S}_i^r(t_j, t^j)$ and $w_j^x \in \mathbb{R}$ is in the interval $[x + h_j, x]$. Define

$$z_j = \sup_{u \in [t_j, t^j)} \frac{\Lambda_i(\bar{Z}^r(u))}{\bar{Z}_i^r(u)} \tag{A.11}$$

and let $\tilde{h}_j = -z_j l/n$. Combine terms and bound the integrand to obtain

$$\left| b_n^r - \sum_{j=0}^{n-1} \langle f' \tilde{h}_j, \bar{\mathcal{Z}}_i^r(t_j) \rangle \right|$$

$$= \left| \sum_{j=0}^{n-1} \langle f(\cdot - \bar{S}_i^r(t_j, t^j)) - f(\cdot) - f'(\cdot)\tilde{h}_j, \bar{\mathcal{Z}}_i^r(t_j) \rangle \right|$$

$$\leq \sum_{j=0}^{n-1} \sup_{x \in \mathbb{R}} |f(x - \bar{S}_i^r(t_j, t^j)) - f(x) - f'(x)\tilde{h}_j| \langle 1, \bar{\mathcal{Z}}_i^r(t_j) \rangle.$$

Apply (A.5) and (A.10) to get

$$\left| b_n^r - \sum_{j=0}^{n-1} \langle f' \tilde{h}_j, \bar{\mathcal{Z}}_i^r(t_j) \rangle \right|$$

$$\leq \sum_{j=0}^{n-1} \sup_{x \in \mathbb{R}} |f'(w_j^x)h_j - f'(x)\tilde{h}_j| \langle 1, \bar{\mathcal{Z}}_i^r(t_j) \rangle \tag{A.12}$$

$$\leq M \sum_{j=0}^{n-1} \sup_{x \in \mathbb{R}} (|f'(w_j^x) - f'(x)||h_j| + |f'(x)||h_j - \tilde{h}_j|).$$

Since $w_j^x \in [x + h_j, x]$ for each $j \in \{0, \ldots, n-1\}$ and $x \in \mathbb{R}$, deduce from the definition of $h_j$, $\bar{S}_i^r(t_j, t^j)$, and from (2.1) and (A.5), that

$$|w_j^x - x| \leq |h_j| = \int_{t_j}^{t^j} \frac{\Lambda_i(\bar{Z}^r(u))}{\bar{Z}_i^r(u)} \, du \leq \frac{l \|C\|}{n\varepsilon}. \tag{A.13}$$

Since $f'$ has compact support, it is uniformly continuous. Hence, there exists a continuous nondecreasing function $\psi_f : \mathbb{R}_+ \to \mathbb{R}_+$ such that $\psi_f(0) = 0$ and for all $h \in \mathbb{R}_+$,

$$\sup_{x \in \mathbb{R}} |f'(x+h) - f'(x)| \leq \psi_f(|h|). \tag{A.14}$$

We deduce from (A.12)–(A.14) that

$$\left| b_n^r - \sum_{j=0}^{n-1} \langle f' \tilde{h}_j, \bar{\mathcal{Z}}_i^r(t_j) \rangle \right|$$
(A.15)
$$\leq M \left( n\psi_f \left( \frac{l\|C\|}{n\varepsilon} \right) \frac{l\|C\|}{n\varepsilon} + \|f'\|_\infty \sum_{j=0}^{n-1} \left( z_j \frac{l}{n} - \bar{S}_i^r(t_j, t^j) \right) \right).$$

Let $\phi_n(u) = \sum_{j=0}^{n-1} z_j 1_{[t_j, t^j)}(u)$ for each $n \in \mathbb{N}$ and $u \in [t_j, t^j]$. Then

(A.16) $$\sum_{j=0}^{n-1} \left( z_j \frac{l}{n} - \bar{S}_i^r(t_j, t^j) \right) = \int_s^t \phi_n(u)\,du - \int_s^t \frac{\Lambda_i(\bar{Z}^r(u))}{\bar{Z}_i^r(u)}\,du.$$

Observe that $\phi_n(u) \to \Lambda_i(\bar{Z}^r(u))\bar{Z}_i^r(u)^{-1}$ as $n \to \infty$, for all $u \in [s,t]$ at which the latter function is continuous, which is at almost every $u$. So, by the bounded convergence theorem, (A.16) converges to zero as $n \to \infty$. This implies, by definition of $\psi_f$, that (A.15) converges to zero as $n \to \infty$. Note that

$$\sum_{j=0}^{n-1} \langle f'\tilde{h}_j, \bar{\mathcal{Z}}_i^r(t_j) \rangle = -\sum_{j=0}^{n-1} \langle f', \bar{\mathcal{Z}}_i^r(t_j) \rangle z_j \frac{l}{n},$$

and that, as $n \to \infty$,

$$-\sum_{j=0}^{n-1} \langle f', \bar{\mathcal{Z}}_i^r(t_j) \rangle z_j \frac{l}{n} \to -\int_s^t \langle f', \bar{\mathcal{Z}}_i^r(u) \rangle \frac{\Lambda_i(\bar{Z}^r(u))}{\bar{Z}_i^r(u)}\,du,$$

by (A.11) and bounded convergence, since the integrand on the right is also continuous at almost every $u$. Conclude that, as $n \to \infty$,

(A.17) $$b_n^r \to -\int_s^t \langle f', \bar{\mathcal{Z}}_i^r(u) \rangle \frac{\Lambda_i(\bar{Z}^r(u))}{\bar{Z}_i^r(u)}\,du.$$

Combining (A.6), (A.9) and (A.17) yields (5.7). □

Department of Mathematics  
University of Virginia  
Charlottesville, Virginia 22903  
USA  
E-mail: gromoll@virginia.edu

Department of Mathematics  
University of California  
San Diego  
9500 Gilman Drive  
La Jolla, California 92093  
USA  
E-mail: williams@math.ucsd.edu